\documentclass{amsart}
\usepackage{amssymb,latexsym}
\usepackage{amscd,amsthm}

\usepackage[all]{xy}
\usepackage[bookmarksnumbered,colorlinks]{hyperref}
\usepackage{tikz}

\usepackage{tikz-cd}

\newtheorem{thm}{\bf Theorem}[section]
\newtheorem{cor}[thm]{\bf Corollary}

\newtheorem{lem}[thm]{\bf Lemma}
\newtheorem{prop}[thm]{\bf Proposition}
\newtheorem{defn}[thm]{\bf Definition}

\newtheorem{rem}[thm]{\bf Remark}

\newtheorem{exmp}[thm]{\bf Example}
\newtheorem{exmps}[thm]{\bf Examples}

\newtheorem{ques}[thm]{\bf Question}

\DeclareMathOperator{\Ext}{Ext}
\DeclareMathOperator{\Hom}{Hom}
\DeclareMathOperator{\Tor}{Tor}

\newcommand{\N}{\mathbb{N}}
\newcommand{\Z}{\mathbb{Z}}
\newcommand{\QZ}{\mathbb{Q/Z}}

\newcommand{\cqfd}
{\hspace{1cm}
\rule{2mm}{2mm}%
\medbreak%
\par%
}

\def\proof{{\parindent0pt {\bf Proof.\ }}}

\def\pd{{\rm pd}}

\def\id{{\rm id}}

\def\resdim{{\rm resdim}}

\def\F{\mathcal{F}_{\rm C}}
\def\P{\mathcal{P}_{\rm C}}

\def\GC{{\rm G}_C}

\def\GCpd{{\rm G_C\,\mbox{-}\,pd}}
\def\PGCFpd{{\rm PG_CF\,\mbox{-}\,pd}}

\def\proj{\mathcal{P}}
\def\flat{\mathcal{F}}
\def\inj{\mathcal{I}}

\def\PGF{\rm PGF}
\def\GF{\rm \mathcal{GF}}
\def\GP{\rm \mathcal{GP}}
\def\GI{\rm \mathcal{GI}}

\def\X{\mathcal{X}}
\def\Y{\mathcal{Y}}
\def\Z{\mathcal{Z}}
\def\G{\mathcal{G}}
\def\A{\mathcal{A}}
\def\B{\mathcal{B}}
\def\Q{\mathcal{Q}}
\def\R{\mathcal{R}}
\def\W{\mathcal{W}}

\def\XX{\textbf{X}}
\def\YY{\textbf{Y}}

\def\pcd{{\rm \mathcal{P}_C\,\mbox{-}\,pd}}

\def\GF{{\rm \mathcal{G}}{\rm \mathcal{F}}}

\def\GCF{{\rm G}_{\rm C}{\rm F}}
\def\GCP{{\rm G}_{\rm C}{\rm P}}
\def\PGCF{{\rm PG}_{\rm C}{\rm F}}

\def\PGCFD{{\rm PG}_{\rm C}{\rm FD}}

\def\Im{{\rm Im}}
\def\Coker{{\rm Coker}}
\def\Ker{{\rm Ker}}
\def\Add{{\rm Add}}

\def\Ext{{\rm Ext}}
\def\Tor{{\rm Tor}}
\def\Hom{{\rm Hom}}
\def\End{{\rm End}}
\def\Prod{{\rm Prod}}
\def\Cogen{{\rm Cogen}}

\def\sup{{\rm sup}}

\begin{document}

\title{Relative PGF modules and dimensions}

\author{Rachid El Maaouy}
\address{}
\email[Rachid El Maaouy]{elmaaouy.rachid@gmail.com}
%\urladdr{http://pages.ramapo.edu/~jgillesp/}

\date{\today}



\begin{abstract}%% what do you actually prove?%%
	
Inspired in part by recent work of \v{S}aroch and \v{S}\v{t}ov\'{\i}\v{c}ek in the setting of Gorenstein homological algebra, we extend the notion of Foxby-Golod ${\rm G_C}$-dimension of finitely generated modules with respect to a semidualizing module $C$ to arbitrary modules over arbitrary rings, with respect to a module $C$ that is not necessarily semidualizing. We call this dimension ${\rm PG_CF}$ dimension and show that it can serve as an alternative definition of the ${\rm G_C}$-projective dimension introduced by Holm and J\o rgensen. Modules with ${\rm PG_CF}$ dimension zero are called ${\rm PG_CF}$ modules. When the module $C$ is nice enough, we show that the class ${\rm PG_CF}(R)$ of these modules is projectively resolving. This enables us to obtain good homological properties of this new dimension. We also show that ${\rm PG_CF}(R)$ is the left-hand side of a complete hereditary cotorsion pair. This yields, from a homotopical perspective, a hereditary Hovey triple where the cofibrant objects coincide with the ${\rm PG_CF}$ modules and the fibrant objects coincide with the modules in the well-known Bass class $\mathcal{B}_C(R)$. 
\end{abstract}

\maketitle

\section{Introduction}

The notion of G-dimension of finitely generated modules over commutative Noetherian rings was introduced by Auslander and studied further by Auslander and Bridger in \cite{AB69}. The significance of this homological dimension has led many researchers to extend it in various directions. Enochs, Jenda and Torrecillas extended it to arbitrary modules over any ring $R$ through the notions of Gorenstein projective and flat modules (see \cite{EJ00}). On the other hand, Foxby \cite{Fox72}, Golod \cite{Gol85} and Holm-J\o rgensen \cite{HJ06}, considered another extension, replacing the base ring with a semidualizing module. Let $R$ be a commutative noetherian ring and $C$ a semidualizing $R$-module (i.e., an $R$-module satisfying $\Hom_R(C,C)\cong R$ and $\Ext^i_R(C,C)=0$ for all $i \geq 1$). For an $R$-module $M$, consider the following three assertions (which are equivalent when restricted to finitely generated modules; see for instance \cite[Theorem 4.4]{Whi10},\cite[Proposition 2.5]{BEGO23a}, and \cite[Proposition 2.14(b)]{BEGO24}):
	\begin{enumerate}
		\item[(1)]  $M$ is totally $C$-reflexive, i.e.,  $M$ is finitely generated and satisfies the following assertions. 
		\begin{enumerate}
			\item[$(a)$] $\Ext_R^{i\geq 1}(M,C)=\Ext_R^{i\geq 1}(\Hom(M,C),C)=0$.
			\item[$(b)$] The natural morphism $M\to \Hom_R(\Hom_R(M,C),C)$ is an isomorphism.
		\end{enumerate}
        \item[(2)]  $M$ is $\GC$-projective, i.e., $M\cong Z_{-1}(\textbf{P})$ is a cycle of an exact complex of  $R$-modules $$\textbf{P}:\cdots \rightarrow P_1\rightarrow P_0 \rightarrow C\otimes_R P_{-1}\rightarrow C\otimes_R P_{-2}\rightarrow\cdots$$ 
        with each $P_i$ projective and  $\Hom_R(\XX,C\otimes_R P)$ is exact for any projective module $P$.
        \item[(3)] $M$ is $\GC$-flat, i.e., $M\cong Z_{-1}(\textbf{F})$ is a cycle of an exact complex of $R$-modules $$\textbf{F}:\cdots \rightarrow F_1\rightarrow F_0 \rightarrow C\otimes_R F_{-1}\rightarrow C\otimes_RF_{-2}\rightarrow\cdots$$ 
        with each $F_i$ flat and $\Hom_R(C,E)\otimes_R \textbf{F}$ is exact for any injective module $E$.
	\end{enumerate}
 
   Foxby \cite[pg. 277]{Fox72} and Golod \cite{Gol85} independently adapted condition (1) to define the $\GC$-dimension using resolutions with respect to totally $C$-reflexive modules. Their $\GC$-dimension shares many homological properties with Auslander's G-dimension. For example, Golod showed that the $\GC$-dimension is a refinement of the projective dimension, i.e., ${\rm G_C-dim}(M)\leq \pd_R(M)$ with equality when $\pd_R(M) <\infty$.  Furthermore, Geng  \cite[Theorem 4.4]{Gen13} extended the Auslander-Bridger formula \cite[4.13(b)]{AB69} to the relative setting. For a comprehensive treatment of the $\GC$-dimension, see Sather-Wagstaff's survey \cite{SW10}. 
 
 On the other hand, Holm and J\o rgensen used condition (2) (resp., (3)) to define the $\GC$-projective (resp., $\GC$-flat) dimension for modules not necessarily finitely generated, extending all previous dimensions. These two dimensions, along with the dual (the so-called $\GC$-injective dimension), form the basis of relative Gorenstein homological algebra with respect to a semidualizing module $C$.
 
In this paper, we consider condition (3) with the modules $F_i$ being projective instead of flat. This gives rise to a new class of modules, which we call $\PGCF$ modules (see Definition \ref{PGCF def}). The main goal of this paper is to show that $\PGCF$ modules behave better both homologically and homotopically than $\GC$-projective modules and can serve as an alternative definition. Our work is heavily influenced by recent work of  \v{S}aroch, and  \v{S}\v{t}ov\'{\i}\v{c}ek in the absolute setting \cite{SS20}. When $C=R$, the $\PGCF$ modules coincide with the PGF modules recently introduced by them. In fact, their work is a major motivation for introducing the $\PGCF$ modules.
\bigskip

Semidualizing modules were first introduced by Foxby \cite{Fox72} under the name PG-modules of rank one over commutative Noetherian rings. Later,  Holm and White defined and studied them over arbitrary associative rings \cite{HW07}. In some situations, requiring a module to be semidualizing is more restrictive than necessary. Our development in this paper is based on modules $C$ with significantly less restrictive conditions than semidualizing modules. This approach was initiated by Bennis, Garc\'ia Rozas, and Oyonarte \cite{BGO16} and has been extensively investigated in the last years. Modules satisfying these conditions are called w-tilting in the case of $\GC$-projectivity  \cite{BGO16} and $w^+$-tilting in the case of $\GC$-flatness \cite{BEGO22}. In our situation,  we call these new modules generalized Wakamatsu tilting, or g-tilting for short (see Definition \ref{g-tilting}). We will show that these modules properly generalize tilting modules in the sense of Bazzoni \cite{Baz04}, tilting modules in the sense of Wakamatsu \cite{Wak04} and semidualizing bimodules in the sense of Holm and White \cite{HW07} (see Propositions \ref{til is g-til} and \ref{semid is g-tilting}). Additionally, the g-tilting modules are particular instances of w-tilting and $w^+$-tilting modules (Proposition \ref{w-til is g-til is w+-til}). Thus, they provide a common ground for studying $\PGCF$, $\GC$-projective and $\GC$-flat modules and dimensions. Our definition of a g-tilting module comes from the nice homological properties of the class $\PGCF(R)$ of all $\PGCF$ modules, as shown in the main results of this paper, which we now describe in more detail.

\bigskip

One of the main properties of building a nice homological dimension over a class of objects in an abelian category with enough projectives is the fact that it is projectively resolving, that is, it contains all projective modules and is closed under extensions and kernels of epimorphisms. In our first main result, we establish (Section \ref{Sec. 4}) this basic homological property for the class of all $\PGCF$ modules.

\bigskip

\noindent {\bf Theorem A.} Let $C$ be a pure projective $R$-module. Then, $C$ is $g$-tilting if and only if $\PGCF(R)$ is projectively resolving and contains the class $\Add_R(C)$. 

\bigskip

Here, $\Add_R(C)$ denotes the class of all $R$-modules that are isomorphic to direct summands of direct sums of copies of $C$. %When $C$ is finitely presented, it coincides with the class of $C$-projective modules, that is, modules of the form $C\otimes_S P$  with $P$ projective left $S$-modules (Lemma \ref{Fc-Pc-Ic}). The main challenge one may face in showing it is in showing closure under extensions.  The $\PGCF$ modules are defined via the tensor functor, and hence the closure under extensions does not follow easily as in the case of objects that are defined with the Hom functor. Our key result to overcome this difficulty is an interesting and recent result by Gillespie and Iacob \cite[Theorem 6]{GI22} that allows us to pass from the tensor functor to the Hom functor (see Proposition \ref{new charac of PGcF}), making this problem easier to solve.

      \bigskip

 As discussed above, the $\GC$-dimension is a refinement of the projective dimension. Motivated by this fact, White asked in \cite[Question 2.15]{Whi10} whether this is the case for the $\GC$-projectivity, that is, whether the $\GC$-projective dimension is a refinement of the projective dimension. It turns out that this is not true in general (see Example \ref{PGCF dim isn't a ref}). So, the question becomes:
  
\begin{ques} When is the $\GC$-projective dimension a refinement of the projective dimension? 
\end{ques}

In our second goal of this paper (Section \ref{Sec. 5}), we study basic homological properties of $\PGCF$ dimension of modules. Then, we show that the $\PGCF$ dimension is a refinement of the $\P$-projective dimension (resolution dimension with respect to the class $\Add_R(C)$), while the $\GC$-projective dimension is a refinement of the $\PGCF$ dimension (Proposition \ref{refinment}). In answer to the above question, we characterize when the $\GC$-projective dimension is a refinement of the projective dimension. For example, this happens if and only if the $\PGCF$ dimension is also a refinement of the projective dimension (Theorem \ref{PGCF dim is a ref}).

\bigskip

From the perspective of relative homological algebra, it is natural to wonder when the class $\PGCF(R)$ provides for (minimal) approximations, i.e., (pre)covers and (pre)envelopes. Since complete cotorsion pairs are one of the main tools in producing such approximations, one may also wonder when $\PGCF(R)$ is part of a complete cotorsion pair. Building on recent ideas of \v{S}aroch, and  \v{S}\v{t}ov\'{\i}\v{c}ek \cite{SS20}, we show,  under the assumption that $C$ is of type ${\rm FP_2}$ (i.e., $C$ and its first syzygy are finitely presented) the following result (Theorem \ref{PGCF-cot-pair} and Corollary \ref{PGCF is covering}). \\

\noindent {\bf Theorem B.} Let $C$ be an $R$-module of type ${\rm FP_2}$. Then, $C$ is g-titling if and only if $(\PGCF(R),\PGCF(R)^\perp)$ is a complete hereditary cotorsion pair cogenerated by a set $\mathcal{S}\cup\{C\}$. Consequently, the following two assertions hold.
\begin{enumerate}
    \item[(a)] $\PGCF(R)$ is a special precovering.
    \item[(b)] $\PGCF(R)$ is covering if and only if  $S:=\End_R(C)$ is left perfect.
\end{enumerate}
    
\bigskip

The investigation from the perspective of homotopy theory connects the class of $\PGCF$ modules with another class of interest: the well-known Bass class $\B_C(R)$. When $C$ is semidualizing, it has been shown that the Bass class, and dually, the Auslander class $\A_C(R)$, are one of the central concepts of (relative) Gorenstein homological algebra. It turns out that the Bass class is also of great importance when it comes to the theory of abelian model structures/Hovey triples. Recall that a triple $(\Q,\W,\R)$ of classes of modules is called a (hereditary) Hovey triple if $(\Q,\W\cap\R)$ and $(\Q\cap\W,\R)$ are (hereditary) complete cotorsion pairs and the class $\W$ is thick (see Section \ref{Sec. 2.4} for details).\\

\noindent {\bf Theorem C.}  Let $C$ be an $R$-module of type ${\rm FP_2}$. Then, $C$ is g-titling if and only if there exists a unique thick class $\W_C(R)$ such that $\left( \PGCF(R),\W_C(R),\B_C(R)\right)$ is a hereditary Hovey triple.  
\bigskip

In this case, we explicitly describe the corresponding abelian model structure in Theorem \ref{PGCF-Hov-trip}, which is cofibrantly generated in the sense of \cite[Sec 2.1.3]{Hov99}.  Unfortunately, unlike the absolute case \cite[Theorem 4.9]{SS20}, we lack an explicit description of the class of weak equivalences/trivial objects. We conclude the paper with two things. First, we give a partial answer to this question (Proposition \ref{trivial objects}). Then, we briefly discuss the dual results of Theorems B and C. 

\section{Preliminaries}

Throughout this paper, $R$ will be an associative (non-necessarily commutative) ring with identity, and all modules will be, unless otherwise specified, unital left $R$-modules. When right $R$-modules need to be used, they will be denoted as $M_R$, while in these cases left $R$-modules will be denoted by $_R M$. Sometimes, left (resp., right) $R$-modules will be identified with right (resp., left) modules over the opposite ring $R^{op}$. The category of all left $R$-modules (resp., , right $R$-modules) will be denoted by $R$-Mod (resp., Mod-$R$). We use $\inj(R)$, $\proj(R)$, and $\flat(R)$ to denote the classes of injective, projective, and flat $R$-modules, respectively. From now on $C$ will stand for an $R$-module, $S$ for its endomorphism ring, $S=\End_R(C)$, and $M^+$ for the character module of an $R$-module $M$, $M^+=\Hom_\mathbb{Z}(M,\QZ)$.

By a subcategory of modules, we will
always mean a full subcategory that is closed under isomorphisms. Any class of modules
will be thought as a (full) subcategory of the category of modules.

In what follows, we shall consider three classes of modules $\mathcal{X},\mathcal{Z}\subseteq R$-Mod and $\mathcal{Y}\subseteq \text{Mod-}R$.

\subsection{(Co)resolutions.} An $\mathcal{X}$-resolution of an $R$-module $M$ is an exact sequence $$\cdots\to X_1\to X_0\to M\to 0$$ where
$X_i\in\mathcal{X}$. An $\mathcal{X}$-coresolution of $M$ is defined dually.

A complex $\textbf{X}$ in
$R$-Mod is called $\left(\mathcal{Y}\otimes_R-\right)$-exact (resp., $\Hom_R(\mathcal{X},-)$-exact,
$ \Hom_R(-,\mathcal{X})$-exact) if $Y\otimes_R\textbf{X}$  (resp., $\Hom_R(X,\textbf{X})$, $\Hom_R(\textbf{X},X)$)
is an exact complex for every $Y\in\mathcal{Y}$ (resp., $X\in\mathcal{X}$).

Consider the following classes:
$$\GF_\Y(\X,-):=\{M| \text{ $M$ has a $(\Y\otimes_R-)$-exact $\X$-resolution} \},$$
$$\GF_\Y(-,\X):=\{M| \text{ $M$ has a $(\Y\otimes_R-)$-exact $\X$-coresolution} \},$$
$$\GP_\Z(\X,-):=\{M| \text{ $M$ has a $\Hom_R(-,\Z)$-exact $\X$-resolution} \},$$
$$\GP_\Z(-,\X):=\{M| \text{ $M$ has a $\Hom_R(-,\Z)$-exact $\X$-coresolution} \},$$
$$\GI_\Z(\X,-):=\{M| \text{ $M$ has a $\Hom_R(\Z,-)$-exact $\X$-resolution} \},$$
$$\GI_\Z(-,\X):=\{M| \text{ $M$ has a $\Hom_R(\Z,-)$-exact $\X$-coresolution} \}.$$
In particular, we set $\mathcal{GF}_\Y(\X)=\mathcal{GF}_\Y(\X,\X),\,\,\,\mathcal{GP}_\Z(\X)=\mathcal{GP}_\Z(\X,\X),\,\,\,\mathcal{GI}_\Z(\X)=\mathcal{GI}_\Z(\X,\X).$

\subsection{Relative (Gorenstein) objects.} For an $R$-module $C$, we use $\Add_R(C)$ (resp., $\Prod_R(C)$) to denote the class of all $R$-modules that are isomorphic to direct summands of direct sums (resp., products) of copies of $C$. We also use $\F(R)$ to denote the class of all modules $M$ such that $M^+\in\Prod_R(C^+)$. Such modules were introduced in \cite{BEGO22} and are called $\F$-flat modules. Following the same terminology, we may call modules in $\Add_R(C)$ (resp., $\Prod_R(C)$) $\P$-projective (resp., $\mathcal{I}_C$-injective).

Recall that an $R$-module $M$ is self-small if the canonical morphism
$\Hom_R(M, M^{(I)}) \to \Hom_R(M,M)^{(I)}$ is an isomorphism for any set $I$. Dually, an $R$-module $M$ is self-cosmall \cite[Definition 2.16]{BEGO24} if the canonical morphisms $(M^+)^I\otimes_R M\to (M^+\otimes_R M)^I \text{ and } M^+\otimes_RM\to \Hom_R(M,M)^+$ are isomorphisms for any set $I$. Any finitely generated module is self-small \cite[Chap I, Exercice 2]{EJ00} and any finitely presented module is self-cosmall \cite[Theorems 3.2.11 and 3.2.22]{EJ00}.
\begin{lem}\label{Fc-Pc-Ic} Let $C$ be an $R$-module.
	\begin{enumerate}
		\item \cite[Proposition 3.1]{BDGO21}  If $C$ is self-small, then $$\Add_R(C)=C\otimes_S \proj(S):=\{C\otimes_S P\,|\, P\in\proj(S)\}.$$
		\item \cite[Corollary 2.19(2)]{BEGO24} If $C$ is self-cosmall, then  $$\Prod_R(C^+)=\Hom_S(C,\inj(S^{op})):=\{\Hom_S(C,E)\,|\, E\in\inj(S^{op})\}.$$
		\item  \cite[Proposition 3.3(2)]{BEGO22} If $C$ is finitely presented, then $$\F(R)=C\otimes_S\flat(S):=\{C\otimes_S F\,|\, F\in\flat(S)\}.$$
	\end{enumerate}
\end{lem}

\begin{defn}[\cite{BGO16,BEGO22}] Let $C$ be an $R$-module. An $R$-module $M$ is said to be $\GC$-projective if there exists a $\Hom_R(-,\Add_R(C))$-exact exact sequence $$\textbf{P}:\cdots \rightarrow P_1\rightarrow P_0 \rightarrow C_{-1}\rightarrow C_{-2}\rightarrow\cdots$$ with each $C_i\in\Add_R(C)$ and $P_j\in\proj(R)$, such that $M\cong \Im(P_0\to C_{-1})$.
	
 An $R$-module $M$ is said to be $\GC$-flat if there exists a  $\left(\Prod_R(C^+)\otimes_R-\right)$-exact  exact sequence $$\textbf{F}:\cdots \rightarrow F_1\rightarrow F_0 \rightarrow C_{-1}\rightarrow C_{-2}\rightarrow\cdots$$ with each $C_i\in\F(R)$ and $F_j\in\flat(R)$, such that $M\cong \Im(F_0\to C_{-1})$.

\end{defn}

The classes of $\GC$-projective and $\GC$-flat $R$-modules are denoted by $\GCP(R)$ and $\GCF(R)$, respectively. It is clear that   $\GCP(R):=\GP_{\Add_R(C)}(\proj(R),\Add_R(C))$ and $\GCF(R):=\GF_{\Prod_R(C^+)}(\flat(R),\F(R)).$

\begin{rem}\item
\begin{enumerate} 
    \item If $C=R$, or more generally, if $C$ is a projective generator, then $\Add_R(C)=\proj(R),$ $\F(R)=\flat(R)$ and $\Prod_R(C^+)=\inj(R^{op}).$ 
In this case,  $\GCP(R)=\GP(R)$ and $\GCF(R)=\GF(R)$ with $\GP(R)$ and $\GF(R)$ denote the classes of the (absolute) Gorenstein projective and flat modules, respectively.

 \item It follows by Lemma \ref{Fc-Pc-Ic} that the $\GC$-projective and $\GC$-flat modules coincide, respectively, with the $C$-Gorenstein projective and $C$-Gorenstein flat modules introduced by Holm and J\o rgensen in \cite{HJ06}.
\end{enumerate}
\end{rem}

 The classes of $\GC$-projective and $\GC$-flat modules behave well when the module $C$ is semidualzing, or more generally when it is w-tilting for the case of $\GC$-projectivity or $w^+$-tilting for the case of $\GC$-flatness.

\begin{defn}[\cite{HW07}, Definition 2.1] An $(R,S)$-bimodule $C$ is semidualizing if:
	\begin{enumerate}
		\item $_RC$ and $C_S$ both admit a degreewise finite projective resolution.
		\item $\Ext_R^{\geq 1}(C,C)=\Ext_S^{\geq 1}(C,C)=0.$
		\item The natural homothety maps $R\rightarrow\Hom_S(C,C)$ and $S \rightarrow \Hom_R(C,C)$  are isomorphisms.
	\end{enumerate}
\end{defn}

By a degreewise finite projective resolution we mean a projective resolution in which every projective module is finitely generated.

As shown in \cite{BGO16} and \cite{BEGO22}, the following two notions properly generalize semdidualizing bimodules.

\begin{defn}  Let $C$ be an $R$-module.
	\begin{enumerate}
\item \cite[Definition 2.1]{BGO16} $C$ is w-tilting if it satisfies the following two properties:
\begin{enumerate}
	\item $C$ is $\Sigma$-self-orthogonal, that is, $\Ext_R^i(C,C^{(I)})=0$ for every set $I$.
	\item $There exists$ a  $\Hom_R(-,\Add_R(C))$-exact $\Add_R(C)$-coresolution $$\textbf{X}:0\rightarrow R\rightarrow C_{-1}\rightarrow C_{-2}\rightarrow\cdots.$$
\end{enumerate}
\item \cite[Definition 4.1]{BEGO22} $C$ is $w^+$-tilting if it satisfies the following two properties:
\begin{enumerate}
	\item $C$ is $\prod$-$\Tor$-orthogonal, that is, $\Tor^R_{i\geq 1}((C^+)^I,C)=0$ for every set $I$.
	\item There exists a  $\left( \Prod_R(C^+)\otimes_R-\right) $-exact $\F(R)$-coresolution $$\textbf{X}:0\rightarrow R\rightarrow C_{-1}\rightarrow C_{-2}\rightarrow\cdots$$
\end{enumerate}
 
	\end{enumerate}
\end{defn}

\bigskip

Associated to an $(R,S)$-bimodule $C$, we have the Foxby classes defined as follows:

$\bullet$ The Bass class $\B_C(R)$ consists of all left $R$-modules $N$ satisfying: 
\begin{enumerate}
    \item $\Ext_R^{\geq 1}(C,N)=\Tor^S_{\geq 1}(C,\Hom_R(C,N))=0,$ and
     \item the canonical map $\nu_N:C\otimes_S\Hom_R(C,N)\to N$ is an isomorphism of $R$-modules.
\end{enumerate}

$\bullet$ The Auslander class $\A_C(S)$ is the class of all left $S$-modules $M$ satisfying:
\begin{enumerate}
    \item $\Tor^S_{\geq 1}(C,M)=\Ext_R^{\geq 1}(C,C\otimes_S M)=0,$ and 
    \item the canonical map $\mu_M:M\to \Hom_R(C,C\otimes_SM)$ is an isomorphism of left $S$-modules.
\end{enumerate}

The Bass and Auslander classes, $\B_C(S)$ and $\A_C (R)$, of right $S$-modules and right $R$-modules can be defined similarly.

\subsection{Cotorsion pairs and approximations.} Given an integer $n\geq 1$, to any given class of $R$-modules $\X$, we associate its right and left $n$-th Ext-orthogonal classes
$$\mathcal{X}^{\perp_n}=\{M\in R\text{-Mod}|\;\Ext^i_R(X,M)=0, \forall X\in\mathcal{X},\forall i=1,..,n \},$$
$$^{\perp_n}\mathcal{X}=\{M\in R\text{-Mod}|\;\Ext^i_R(M,X)=0, \forall X\in\mathcal{X},\forall i=1,..,n \}.$$
In particular, we set $\mathcal{X}^{\perp}=\mathcal{X}^{\perp_1},$   $\mathcal{X}^{\perp_\infty}=\cap_{n\geq 1}\mathcal{X}^{\perp_n},$ 
$^{\perp}\mathcal{X}=\;^{\perp_1}\mathcal{X}$
and $ ^{\perp_\infty}\mathcal{X}=\cap_{n\geq 1}\;^{\perp_n}\mathcal{X}$
%$$\xymatrix{  \mathcal{X}^{\perp}=\mathcal{X}^{\perp_1},  &  \mathcal{X}^{\perp_\infty}=\cap_{n\geq 1}\mathcal{X}^{\perp_n},  &  ^{\perp}\mathcal{X}=\;^{\perp_1}\mathcal{X},     &   ^{\perp_\infty}\mathcal{X}=\cap_{n\geq 1}\;^{\perp_n}\mathcal{X}}$$
and if $\X=\{X\}$, we simply write $\X^{\perp_n}=X^{\perp_n}$ and $^{\perp_n}\X=\;^{\perp_n}X.$

An $\mathcal{X}$-precover of a module $M$ is a morphism $f:X\to M$ with $X\in \mathcal{X}$, in such a way that the map $f_*:\Hom_R(X',X)\to \Hom_R(X',M)$ is surjective for every $X'\in \mathcal{X}$. An $\mathcal{X}$-precover is called an $\mathcal{X}$-cover if every endomorphism $g:X\to X$ such that $fg=f$ is an automorphism of $X$. If every module has an $\mathcal{X}$-(pre)cover, then the class $\mathcal{X}$ is said to be (pre)covering.  An $\mathcal{X}$-precover is called special if it is epimorphism and $\Ker f\in\X^\perp$. $\mathcal{X}$-(special pre)envelopes can be defined dually.

An important approach for producing (pre)covers and (pre)envelopes is via complete cotorsion pairs. Recall that a pair $\left(\mathcal{A},\mathcal{B}\right)$ of classes of modules is called a cotorsion pair if $\mathcal{A}^\perp=\mathcal{B}$ and $\mathcal{A}=\,^\perp\mathcal{B}$. A cotorsion pair $\left(\mathcal{A},\mathcal{B}\right)$ is said to be:
\begin{itemize}
    \item Hereditary if $\Ext^i_R(X,Y)=0$ for every $X\in\X$, $Y\in\Y$ and $i\geq 1$. %$\A$ is projectively resolving, or equivalently, if $\mathcal{B}$ is injectively coresolving. 
\item  Complete if any $R$-module has a special $\B$-preenvelope, or equivalently, any
$R$-module has a special $\A$-precover. 
\item Perfect if every module has an $\mathcal{A}$-cover and a $\mathcal{B}$-envelope. 
\item Cogenerated by a set if there is a set $\mathcal{S}$ such that $\mathcal{B}=\mathcal{S}^\perp$.
\end{itemize}

It is well known that a perfect cotorsion pair $(\mathcal{A},\mathcal{B})$ is complete.  The converse holds when $\mathcal{A}$ is closed under direct limits \cite[Corollary 5.32]{GT12}. We also know by \cite[Theorem 6.11]{GT12} that any cotorsion pair cogenerated by a set is complete.  Standard references for approximations and cotorsion pairs include \cite{EJ00} and \cite{GT12}. 

\subsection{Abelian model structures.} \label{Sec. 2.4} Recall that an abelian model structure is a Quillen model structure, that is, three classes of morphisms of $\A$ called cofibrations, fibrations and weak equivalences, satisfying a set of axioms (\cite[Definition 1.1.3.]{Hov99}), and the following two assertions:
\begin{enumerate}
\item[(a)] A morphism is a (trivial) cofibration if and only if it is a monomorphism with (trivially) cofibrant cokernel. 
\item [(b)] A morphism is a (trivial) fibration if and only if it is an epimorphism with (trivially) fibrant kernel.
\end{enumerate}

Hovey showed in \cite[Theorem 2.2]{Hov02} that an abelian model structure on $R$-Mod (in fact, on any abelian category), is equivalent to a triple $\mathcal{M}=(\mathcal{Q},\mathcal{W},\mathcal{R})$ of classes of modules in $R$-Mod such that $\mathcal{W}$ is thick and $(\mathcal{Q},\mathcal{W}\cap\mathcal{R})$ and $(\mathcal{Q}\cap\mathcal{W},\mathcal{R})$ are complete cotorsion pairs. In this case, $\Q$ is precisely the class of cofibrant objects, $\R$ is precisely the class of fibrant objects, and $\W$ is the class of trivial objects of the model structure.

\section{Homological closure properties.}\label{Sec. 3}
For the sake of better readability, in this section we prove and collect some preliminary results that are of independent interest. To build a good homological dimension with respect to a given class, this class must satisfy some basic homological properties. We discuss in this section these properties for the classes $\GF_\Y(\X,-)$, $\GF_\Y(-,\X)$, $\GP_\Z(-,\X)$ and $\GP_\Z(\X,-)$. For the two classes $\GI_\Z(\X,-)$ and $\GI_\Z(-,\X)$, the corresponding results and their proofs are dual.

\begin{lem}(\textbf{Closure under extensions}) \label{cl-ext} Assume that $\mathcal{X}$ is closed under finite direct sums and let $\mathcal{E}:0\to M'\to M\to M''\to 0$ be an exact sequence of $R$-modules.
\begin{enumerate} 
	\item Assume that $\Hom_R(\mathcal{X},\mathcal{E})$ is exact.  The following two assertions hold:
	\begin{enumerate}
\item[$(i)$] If $M',M''\in\GF_\Y(\X,-)$, then $M\in\GF_\Y(\X,-)$. 

\item[$(ii)$] If $M',M''\in \GP_\Z(\X,-)$, then $M\in \GP_\Z(\X,-)$.
    \end{enumerate}
	\item   Assume that $\Hom_R(\mathcal{E},\mathcal{X})$ is exact. The following two assertions hold: 
		\begin{enumerate}
    \item[$(i)$] If $M',M''\in\GF_\Y(-,\X)$, then $M\in\GF_\Y(-,\X)$. 

	\item[$(ii)$] If $M',M''\in \GP_\Z(-,\X)$, then $M\in \GP_\Z(-,\X)$.

	\end{enumerate} 
\end{enumerate}
\end{lem}
\proof  We only deal with (1)(i) since $(2)(ii)$  and $(2)$ can be shown by a similar argument. 

Assume that there exist two $(\Y\otimes_R-)$-exact $\mathcal{X}$-resolutions
$$\XX': \cdots \to X_1'\to  X_0'\to M' \to 0 \,\,\text{ and } \,\,\XX'':  \cdots \to X_1''\to  X_0''\to M'' \to 0.$$
Let $K_i'=\Ker(X'_i\to X'_{i-1})$ and $K_i''=\Ker(X''_i\to X''_{i-1})$ with $X'_{-1}:=M'$ and $X''_{-1}:=M''$. By the Horseshoe Lemma, it is standard to construct the following commutative diagram with exact rows and columns $$\xymatrix{  &  0 \ar[d]    & 0 \ar[d]  &  0\ar[d]  &  \\
	0 \ar[r] &  K_0' \ar[r] \ar[d]    & K_0 \ar[d] \ar[r] &  K''_0 \ar[d] \ar[r] & 0\\ 
	0 \ar[r] & X'_0 \ar[r] \ar[d]    & X_0'\oplus X_0'' \ar[d] \ar[r] &  X''_1  \ar[d] \ar[r] & 0 \\  
	0 \ar[r] & M' \ar[r] \ar[d]    & M \ar[r] \ar[d] &  M'' \ar[r] \ar[d] & 0  \\ 
	&  0    & 0   &  0.   &}$$
Let $Y\in\Y$. If we apply the functor $Y\otimes_R-$ to this diagram and using the Snake Lemma, we get a $(Y\otimes_R-)$-exact short exact sequence $$0\to K_0\to X_0'\oplus X_0''\to M\to 0.$$ 
Note that $\mathcal{E}_0:\, 0\to K'_0\to K_0\to K''_0\to 0$ is $\Hom_R(\X,-)$-exact. So,
 repeating the same argument infinitely many times, we obtain a short exact sequence of exact complexes 
 $$ 0\to \XX'\to \XX\to \XX''\to 0$$
with a $(\Y\otimes_R-)$-exact $\X$-resolution of $M$:
$$\XX: \cdots \to X_1'\oplus X_1''\to  X_0'\oplus  X_0''\to M \to 0.$$ 
Thus, $M\in \GF_\Y(\X,-)$. \cqfd

\begin{defn}Assume that  $\X,\Z\subseteq\mathcal{G}$ are three classes of $R$-modules.
	
	\begin{enumerate}
		\item $\mathcal{G}$ is said to be \textbf{left $\X$-periodic} if, for any $R$-module $M$, $M \in \mathcal{G}$ if and only if $\X$ is a generator for $\G$, i.e., there exists a short exact sequence of $R$-modules $0\to G\to X\to M\to 0$ with $X\in \X$ and $G\in \mathcal{G}$.

		\item $\mathcal{G}$ is said to be \textbf{right $\X$-periodic} if, for any $R$-module $M$, $M\in \mathcal{G}$ if and only if $\X$ is a cogenerator for $\G$, i.e., there exists a short exact sequence of $R$-modules $0\to M\to X\to G\to 0$ with $X\in \X$ and $G\in \mathcal{G}$
		\item $\mathcal{G}$ is said to be \textbf{$(\X,\Z)$-periodic} if it is left $\X$-periodic and right $\Z$-periodic.
		
		When $\Z=\X$, we simply say that $\G$ is \textbf{$\X$-periodic}.
	\end{enumerate}
\end{defn}

Recall \cite{Wan20} that $\X$ is a $\Y$-pure generator (resp., cogenerator) for $\G$ if any $R$-module $M$ fits into a $(\Y\otimes-)$-exact exact sequence as in (1) (resp., (2)).

\begin{exmps}(Gorenstein) homological algebra is a rich source of examples of periodic classes. For instance, we have the following examples:
\begin{enumerate}
\item By definition, the classes of Gorenstein projective, Gorenstein injective and Gorenstein flat modules are $\proj(R)$-periodic, $\inj(R)$-periodic and $\flat(R)$-periodic classes, respectively.

\item In the relative setting, the class of $\GC$-projective modules is $(\proj(R),\Add_R(C))$-periodic when $C$ is w-tilting \cite[Propositions 2.8, 2.9 and 2.11]{BGO16}. Dually, the class of ${\rm G_U}$-injective modules in the sense of \cite[Definition 4.1]{BGO16} is $(\Prod_R(U),\inj(R))$-periodic when $U$ is weakly cotilting (see \cite[Section 4]{BGO16}).
	\end{enumerate}
\end{exmps}

\begin{lem} (\textbf{Closure under epikernels/monocokernels}) \label{cl-epikernels} Let $\mathcal{G}$ be a class of modules closed under extensions and $\X\subseteq \mathcal{G}$.
	\begin{enumerate}		
\item If $\G$ is right $\X$-periodic, then  $\G$ is closed under kernels of epimorphisms. 
\item If $\G$ is left $\X$-periodic, then $\G$ is closed under cokernels of monomorphisms. 
	\end{enumerate}
\end{lem}
\proof 1. Let $0\to M\to N\to L\to 0$ be a short exact sequence of $R$-modules with $N,\ L\in \mathcal{G}$. By hypothesis, there exists a short exact sequence $0\rightarrow N\rightarrow X\rightarrow G\rightarrow 0$ with $X\in\X$ and $G\in \mathcal{G}$. Consider the pushout diagram $$\xymatrix{ & & 0\ar[d] & 0\ar[d] & \\ 0\ar[r] & M\ar@{=}[d] \ar[r] & N\ar[d] \ar[r] & L\ar@{-->}[d] \ar[r] & 0 \\ 0\ar[r]& M\ar[r] & X\ar@{-->}[r] \ar[d]& H\ar[r] \ar[d] & 0 \\ &  &G\ar@{=}[r] \ar[d] & G \ar[d]  & \\ & & 0 & 0 & }$$
As $\mathcal{G}$ is closed under extensions, $H\in\mathcal{G}$ and hence, using hypothesis and the middle row, $M\in\mathcal{G}$. \cqfd

Using Eilenberg’s trick, the following closure follows the same argument used by Holm in \cite[Proposition 1.4]{Hol04}.
\begin{lem} (\textbf{Closure under direct summands}) \label{cl-dir-summnds} Let $\G$ be a class of $R$-modules satisfying one of the following assertions:
\begin{enumerate}
\item[$(a)$] $\G$ is closed under kernels of epimorphism and countable direct sums.
\item[$(b)$] $\G$ is closed under cokernels of monomorphism and countable direct products.
\end{enumerate}

Then, $\G$ is closed under direct summands. 
\end{lem}

\bigskip

Given a class $\mathcal{S}$ of $R$-modules, recall that an $R$-module $M$ is a transfinite extension of $\mathcal{S}$ (or $\mathcal{S}$-filtered) if $M \cong \varinjlim_{\alpha<\lambda}S_\alpha$ for some well-ordered direct system $(S_\alpha|\,\alpha\leq \lambda)$ indexed by an ordinal $\lambda$, such that 
\begin{itemize}
    \item  For every $\alpha+1\leq \lambda$, the morphism $S_\alpha\to S_{\alpha+1}$ is a monomorphism.
    \item  $S_\alpha\cong \varinjlim_{\beta<\alpha}S_\beta$ for any limit ordinal $\alpha\leq \lambda$.
    \item  $S_0=0$ and $S_{\alpha+1}/S_\alpha:=\Coker(S_\alpha\to S_{\alpha+1})\in\mathcal{S}$  for every $\alpha+1\leq \lambda$.
\end{itemize}

A class $\mathcal{S}$ of modules is said to be closed under transfinite extensions provided that every $\mathcal{S}$-filtered module belongs to $\mathcal{S}$.  Clearly, this implies that $\mathcal{S}$ is closed under extensions and finite direct sums.

%-(Similarly, we define continuous direct systems of exact sequences for well-ordered direct systems of short exact sequences of modules)

\begin{lem} (\textbf{Closure under transfinite extensions (1)}) \label{transf-ext 1} Assume that $\mathcal{X}$ is closed under transfinite extensions. 
The following two assertions hold:
	\begin{enumerate} 
		\item If  $\Ext^1_R( \X,\GF_\Y(\X,-))=0$, then $\GF_\Y(\X,-)$ is closed under transfinite extensions.
		\item  If $\Ext^1_R(\GF_\Y(\X,-),\X)=0$, then  $\GF_\Y(-,\X)$ is closed under transfinite extensions.
		\item   If $\Ext^1_R(\X,\GP_\Z(\X,-))=0$, then $\GP_\Z(\X,-)$ is closed under transfinite extensions.
	    \item  If $\Ext^1_R(\GP_\Z(\X,-),\X)=0$, then $\GP_\Z(-,\X)$ is closed under transfinite extensions.
	\end{enumerate}
\end{lem}
\proof $1.$ Let $M$ be an $R$-module and $(M_\alpha | \alpha\leq \kappa)$ be a $\mathcal{GF}_\Y(\X,-)$-filtration of $M$. By induction on $\alpha\leq \kappa$, we will prove that each $M_\alpha\in \mathcal{GF}_\Y(\X,-)$. %That is, for any $\alpha\leq \kappa$ we construct. 

 If $\alpha=0$, we set $\textbf{X}_\alpha=0$. Assume that $\alpha>0$ and we have constructed a $\left( \mathcal{Y}\otimes_R-\right)$-exact $\mathcal{X}$-resolution of $M_\alpha$:
$$\XX_\alpha: \cdots\to X_{\alpha,1}\to X_{\alpha,0}\to M_\alpha\to 0.$$
By hypothesis, there exists a $\left(\mathcal{Y}\otimes_R-\right)$-exact $\mathcal{X}$-resolution of $M_{\alpha+1}/ M_\alpha$:
$$\XX'_\alpha: \cdots\to X'_{\alpha,1}\to X'_{\alpha,0}\to M_{\alpha+1}/ M_\alpha\to 0.$$
Clearly, we have a short exact sequence of $R$-modules
$$\mathcal{E}_\alpha: 0\to M_\alpha\to M_{\alpha+1}\to M_{\alpha+1}/M_{\alpha}\to 0,$$
 which is $\Hom_R(\X,-)$-exact since $\Ext^1_R( \X,M_{\alpha})=0$. So, we obtain (as in Lemma \ref{cl-ext}(1)) a short exact sequence of exact complexes 
 $$ 0\to \XX_{\alpha}\to \XX_{\alpha+1}\to \XX'_{\alpha}\to 0$$
 where
$$\XX_{\alpha+1}: \cdots \to X_{\alpha,1}\oplus X'_{\alpha,1}\to  X_{\alpha,0}\oplus X'_{\alpha,0}\to M_{\alpha+1} \to 0$$ 
is a $(\Y\otimes_R-)$-exact $\X$-resolution of $M_{\alpha+1}$.  Therefore, we can set
$$\begin{cases}
X_{\alpha+1,k}:=X_{\alpha,k}\oplus X'_{\alpha,k} & \text{ for any $\alpha+1\leq \kappa$} \\ 
X_{\alpha,k}:=\varinjlim_{\beta<\alpha} X_{\beta,k} & \text{ for any limit ordinal $\alpha\leq \kappa$}
\end{cases}$$
Since $\X$ is closed under transfinite extensions, each $X_k:=\varinjlim_{\alpha<\kappa}\X_{\alpha,k}\in\X$. Finally, taking the direct limit of all the $\XX_\alpha$, we get an $\mathcal{X}$-resolution of $M$: 
$$\XX: \cdots\to X_1\to X_0\to M \to 0,$$
which is  $\left( \mathcal{Y}\otimes_R-\right)$-exact, since the direct limit is exact and commutes with the tensor products. Therefore, $M\in \mathcal{GF}_\Y(\X,-)$.
\cqfd

\begin{lem}(\textbf{Closure under transfinite extensions (2)})\label{transf-ext 2} Assume that $\X$ is self-orthogonal and closed under direct sums. Then $\X$ is closed under transfinite extensions. %In partuclar, the class  $\Add_R(C)$ is closed under transfinite  extensions for any $\Sigma$-self-ortohogonal module $C$.
\end{lem}
\proof Let $X$ be an $R$-module and $(X_\alpha\;|\; \alpha\leq \kappa)$ be an $\X$-filtration of $X$. We may assume that each $X_\alpha$ is a submodule of $X_{\alpha+1}$. By induction on $\alpha\leq \kappa$, we will prove that   

$$X_\alpha \cong  \bigoplus_{\beta<\alpha}(X_{\beta+1}/X_\beta).$$

If $\alpha=0$, there is nothing to show. So, we assume that $\alpha>0$ and  $X_\alpha\cong  \bigoplus_{\beta<\alpha}(X_{\beta+1}/X_\beta)$. As $\X$ self-ortohogonal and $X_\alpha,X_{\alpha+1}/X_\alpha\in\X$, the short exact sequence of $R$-modules splits 
$$\mathcal{E}_\alpha:\;0\to X_\alpha\stackrel{\subseteq}{\longrightarrow} X_{\alpha+1}\to X_{\alpha+1}/X_\alpha \to 0.$$
Then, 
$X_{\alpha+1}\cong X_\alpha \oplus X_{\alpha+1}/X_\alpha \cong  \bigoplus_{\beta<\alpha}(X_{\beta+1}/X_\beta)\oplus (X_{\alpha+1}/X_\alpha)\cong \bigoplus_{\beta<\alpha+1}(X_{\beta+1}/X_\beta).$

Assume now that $\alpha$ is a limit ordinal. By induction hypothsis, $X_\beta= \bigoplus_{\gamma<\beta} (X_{\gamma+1}/X_\gamma)$ for any $\beta<\alpha.$ Then, for each $\beta<\alpha,$ there exists a submodule  $Y_\beta$ of $X_{\beta+1}$ with $Y_\beta \cong X_{\beta+1}/X_{\beta}$ and $X_\beta= \bigoplus_{\gamma<\beta} Y_\gamma$. Let us show that $$X_\alpha=\sum_{\beta<\alpha} Y_\beta =\bigoplus_{\beta<\alpha} Y_\beta $$
with $\sum_{\beta<\alpha} Y_\beta $ the submodule generated by all $X_\beta$. Since $$X_\alpha=\bigcup_{\beta<\alpha} X_\beta \text{ and } X_\beta= \bigoplus_{\gamma<\beta} Y_\gamma= \sum_{\gamma<\beta} Y_\gamma,  \,\,\,\forall\beta<\alpha,$$  
we get $X_\alpha= \sum_{\beta<\alpha} Y_\beta$. Moreover, to see that the $Y_\beta$ generate their direct sum, suppose that $y_{\beta_1}+\cdots +y_{\beta_n}=0$ is an equation with each $y_{\beta_i}\in Y_{\beta_i}$ and $\beta_1<\cdots <\beta_n<\alpha$. Since $X_{\beta_n+1}= \bigoplus_{\beta\leq \beta_n} Y_\beta$, we get that $y_{\beta_1}=\cdots =y_{\beta_n}=0$. Therefore,  
$$X_\alpha=\bigoplus_{\beta<\alpha} Y_\beta \cong \bigoplus_{\beta<\alpha} (X_{\beta+1}/X_{\beta}).$$
By transfinite induction, this happens for any $\alpha\leq \kappa$, in particular for $\alpha=\kappa$. Thus, $$X=X_\kappa\cong \bigoplus_{\alpha<\kappa} (X_{\alpha+1}/X_{\alpha})\in\X. \text{  \cqfd}$$

Now, we discuss another interesting closure property that is needed in this paper. This closure problem was originally raised by Sather-Wagstaff, Sharif and White in \cite{SSW08}. Namely, they asked: what is the result of the iteration of constructing Gorenstein objects? In other words, if one computes a class $G(\X)$ of Gorenstein modules by taking as the base class, a given class $\X$, then inductively defines $G^n(\X):=G(G^{n-1}(\X))$ for any integer $n\geq 2$, when does this iterative process stabilize in the sense that $G^n(\X)=G(\X)$ for some $n$? 

Many classes of modules have this kind of stability. Recently, Wang provided a general approach to this problem, unifying the stability results for various classes ~\cite{Wan20}.  According to the terminology we have been using so far, we summarize Wang's results needed in this paper as follows.

\begin{lem} [\cite{Wan20}, Lemmas 2.1, 2.2, 2.14 and 2.15] \label{stability problem} Assume that $\X\subseteq \mathcal{G}$ and $\mathcal{Z}$ are three classes of $R$-modules such that   $\mathcal{G}$ is closed under extensions. 
	\begin{enumerate}
		\item If $\X$ is a $\Y$-pure generator for $\G$, then  $\GF_\Y(\X,-)=\GF_\Y(\G,-)$.
		\item If $\X$ is a $\Y$-pure cogenerator for $\G$, then $\mathcal{GF}_\Y(-,\X)=\mathcal{GF}_\Y(-,\G)$.
		\item If $\X$ is a $\Z$-proper generator for $\G$, then $\mathcal{GP}_\Z(\X,-)=\mathcal{GP}_\Z(\G,-)$.
		\item If $\X$ is a $\Z$-proper cogenerator for $\G$, then $\mathcal{GP}_\Z(-,\X)=\mathcal{GP}_\Z(-,\G).$
		%\item  $\mathcal{GI}_\Z(\X,-)=\mathcal{GI}_\Z(\G,-)$
		%\item $\mathcal{GI}_\Z(-,\X)=\mathcal{GI}_\Z(-,\G)$
	\end{enumerate}
\end{lem}
\bigskip

Another key concept needed in the paper is that of definable classes. Recall \cite[Theorem 3.4.7]{Pre09} that a class $\mathcal{D}$ of modules is definable if it is closed under products, direct limits and pure submodules. Given a class of modules $\X$, we shall denote by:

$\bullet$ $\langle\X\rangle$ the definable closure of $\X$, that is, the smallest definable subcategory containing $\X$.

$\bullet$ $\mathcal{X}^p$ will stand for the subcategory of pure submodules of objects in $\mathcal{X}$.

$\bullet$ $\prod\mathcal{X}$ will denote the subcategory of direct products of modules in $\mathcal{X}$.

$\bullet$ $\Prod_R(\mathcal{X})$ will denote the subcategory of direct summands of modules in $\prod\mathcal{X}$.

$\bullet$ ${\rm Cogen}_*(\X)$ the class
of modules cogenerated by $\X$, i.e. the class of all pure submodules of products of modules from $\X$. Note that ${\rm Cogen}_*(\X)=(\prod\mathcal{X})^p$.

In case $\X=\{X\}$, we simply write $\langle\X\rangle=\langle X\rangle$ and $\Cogen^*_R(\mathcal{X})=\Cogen^*_R(X)$.

\begin{lem}\label{definability of Add and Prod} Let $X$ be an $R$-module. Then, $\langle X\rangle=\langle\Add_R(X)\rangle =\langle\Prod_R(X)\rangle=\langle X^{++}\rangle.$
\end{lem}
\proof 1.  The inclusion $\langle X\rangle\subseteq \langle\Add_R(X)\rangle$ is immediate as $X\in \Add_R(X)\subseteq \langle\Add_R(X)\rangle$.  For the other inclusion, we only need to show that $\Add_R(X)\subseteq \langle X\rangle$.  Given a module $Y\in\Add_R(X)$, there exist a module $Z$ and a set $I$ such that $Y\oplus Z=X^{(I)}$. Since $X\in \langle X\rangle$ and $\langle X\rangle$ is closed under direct sums and summands, $Y\in \langle X\rangle$. 

The equality $\langle X\rangle =\langle\Prod_R(X)\rangle $ can be proved by a similar argument.

Finally, it remains show the last equality  $\langle X\rangle=\langle X^{++}\rangle$. Since $\langle X^{++}\rangle$ is closed under pure submodules and the canonical morphism $X\hookrightarrow X^{++}$ is a pure monomorphism, $X\in \langle X^{++}\rangle$. Hence, $\langle X\rangle\subseteq \langle X^{++}\rangle$. Conversely, by \cite[Lemma 1.1]{EIP20}, $X^{++}\in \langle X\rangle$ since $X\in \langle X^{++}\rangle$ and hence $\langle X^{++} \rangle\subseteq  \langle X\rangle$. \cqfd

It seems that the following description of definable closure of a class of modules $\X$ is well known to researchers in model theory. As we were unable to locate a suitable reference, we provide here a proof due to \v{S}aroch (private communication).

\begin{lem}\label{definability} Let $\mathcal{D}$ denote the class of $R$-modules obtained after closing $\X$ under direct products, then under pure submodules and, finally, under pure quotients. Then $\mathcal{D}=\langle \X \rangle$.

\end{lem}
\proof Let $M\in \langle \X \rangle$. By Rothmeler \cite[Fact 7.7]{Rot97} (see also the proof of \cite[Lemma 2.9]{Her14}), $M$ is a pure submodule of a reduced product of modules from $\X$. By \cite[Lemma 3.3.1]{Pre09}, this implies that $M$ is a pure submodule of a pure quotient of a product of modules from $\X$. So, there exists a pure exact sequence $0\to K\to \prod_{i\in I} X_i\to N\to 0$, where $(X_i)_{i\in I}$ is a family of modules in $\X$, such that $M$ is pure in $N$. Now, we form the pull-back of the pure epimorphism from the short exact sequence and the pure embedding of $M$ into $N$. It is straightforward to show that the two
resulting pullback maps are again pure epimorphism and pure monomorphism. Thus, we get $M$ as a pure quotient of a pure submodule of $\prod_{i\in I} X_i$, that is, $M\in\mathcal{D}$.  The other inclusion is clear from the definition and \cite[Theorem 3.4.8]{Pre09}. \cqfd

\begin{lem}\label{def-tensor} Let $\Y$ be a class of right $R$-modules.
	\begin{enumerate}
		\item  $\langle \Y\rangle=\{M| \;0\to F_1\to F_0\to M\to 0 \text{ is a pure s.e.s with } F_i\in \Cogen_*(\Y)\}.$
		
		\item Assume that $\Y$ is closed under products and let $\XX$ be an exact complex of $R$-modules.
		
		\begin{enumerate}
\item[$(a)$] $Y\otimes_R\XX$ is exact for every $Y\in \Y$.
\item[$(b)$] $Y\otimes_R\XX$ is exact for every $Y\in \Cogen_*(\Y)$.
\item[$(c)$] $Y\otimes_R\XX$ is exact for every $Y\in \langle \Y\rangle$.
\end{enumerate}

%Consequently,     $\langle\Y\rangle^\intercal=\Cogen_*(\Y)^\intercal=\Y^\intercal$.

	\end{enumerate}
\end{lem}

\proof (1) The closure under products of $\Y$ is $\prod\Y$. Then the closure under pure
submodules of $\prod\Y$ is $(\prod\Y)^p=\Cogen^*(\Y)$. Finally, closing $\Cogen^*(\Y)$ under pure quotients, we obtain our desired quality by Lemma \ref{definability}.

(2) The implications $(c) \Rightarrow (b) \Rightarrow (a)$ are clear since $\Y\subseteq \Cogen_*(\Y) \subseteq \langle\Y\rangle$.

$(a) \Rightarrow (b)$ Let $Y\in \Cogen_*(\Y)$. Then there exists a pure short exact sequence of right $R$-modules $$0\to Y\to E\to L\to 0$$
with $E\in \prod\Y=\Y$.  Then, $0\to L^+\to E^+\to Y^+\to 0$ splits. Applying the functor $\Hom_R(\XX,-)$ to this split short exact sequence,  we get a split short exact sequence of complexes
$$0\to \Hom_R(\XX,L^+)\to \Hom_R(\XX,E^+)\to \Hom_R(\XX,Y^+)\to 0.$$
By hypothesis, $E\otimes_R\XX$ is exact, so is $\Hom_R(\XX,E^+)\cong (E\otimes_R\XX)^+$. Therefore,  $Y\otimes_R\XX$ is exact since $(Y\otimes_R\XX)^+\cong \Hom_R(\XX,L^+)$ is a direct summand of $\Hom_R(\XX,E^+)$, which is exact.

$(b) \Rightarrow (c)$ Let $Y\in \langle \Y\rangle$. By (1), there exists a pure short exact sequence
$$0\to F_1\to F_0\to Y\to 0$$
with each $F_i\in \Cogen_*(\Y)$. We get after applying the functor $-\otimes_R\YY$, a short exact sequence of complexes
$$0\to F_1\otimes_R\XX\to F_0\otimes_R\XX\to Y\otimes_R\XX\to 0.$$
By hypothesis, each $F_i\otimes_R\YY$ is exact, so is $X\otimes_R\YY$. \cqfd

%Finally, we show the equality $\langle\Y\rangle^\intercal=\Y^\intercal$. The inclusion $\langle\Y\rangle^\intercal\subseteq \Y^\intercal$ is obvious.  Let $X\in \Y^\intercal$. Then, any projective resolution of $X$ is $(\Y\otimes-)$-exact. Therefore, by the above equivalences, any projective resolution of $X$ is $(\langle\Y\rangle\otimes_R-)$-exact. Hence, $\Tor_1^R(\langle\Y\rangle,X)=0$, that is, $X\in \langle\X\rangle^\intercal$. 

\bigskip

Let $\mathcal{A}$ and $\mathcal{B}$ be two abelian categories and $(F,G):\A \rightleftarrows \B$ an adjoint pair. If $\tau:\Hom_\A(-,G(-)) \to \Hom_\B(F(-),-)$ is the natural equivalence given by the adjoint pair, then the unit and the counit are natural transformations.
$$\mu: 1_\A\to GF,\, \mu_A= \tau^{-1}_{A,F(A)}(1_{F(A)}) \text{ and } \nu: FG \to 1_\B,\, \nu_B= \tau_{G(B),B}(1_{G(B)}).$$

Keeping these notation in mind, we conclude this section with the following lemma, which is a categorical generalization of \cite[Theorems 2.2 and 2.3]{Mao19}.

\begin{lem}\label{cov-env transfer} Let $\mathcal{D}\subseteq \A$ and $\mathcal{C}\subseteq \B$ be two classes of objects. Then, the following statements hold for an object $N\in \A$ and $M\in \B$.

\begin{enumerate}
    \item Assume that $\mu_D$ is an isomorphism for any objet $D\in\mathcal{D}$. Then $M$ has a $F(\mathcal{D})$-(pre)cover if and only if $G(M)$ has a $\mathcal{D}$-(pre)cover.
    \item  Assume that $\nu_C$ is an isomorphism for any objet $C\in\mathcal{C}$. Then $N$ has a $G(\mathcal{C})$-(pre)cover if and only if $F(N)$ has a $\mathcal{C}$-(pre)cover.
\end{enumerate}
    
\end{lem}

\proof (2) is dual to (1), while both directions in $(1)$ are similar, so we only show one direction.

(1) ($\Rightarrow$) Let $f:F(D)\to M$ be a $F(\mathcal{D})$-precover. We claim that the composite 
$$g:= G(f)\mu_D: D \stackrel{\mu_D}{\longrightarrow} GF(D) \stackrel{G(f)}{\longrightarrow} 
 G(M)$$
 is a $\mathcal{D}$-precover. Let $D'\in\mathcal{D}$ be an object. Using the adjoint isomorphism and the isomorphism $\mu_D: D\to GF(D)$, we obtain a commutative diagram
$$\xymatrixcolsep{7pc}\xymatrix{  \Hom_\B(F(D'),F(D)) \ar[r]^{\Hom_\B(F(D'),f)} \ar[d]^{\cong}  & \Hom_\B(F(D'),M) \ar[d]^{\cong} \\
\Hom_\A(D',D) \ar[r]^{\Hom_\A(D',G(f))} &  \Hom_\A(D',G(M))}.$$
Since $f$ is an $F(\mathcal{D)}$-precover,  $\Hom_\B(F(D'),f)$ is surjective, and so is $\Hom_\A(D',G(f))$. Hence, $g$ is a precover.

Assume now that $f:F(D)\to M$ be a $F(\mathcal{D})$-cover. Let $h:D\to D$ be an endomorphism such that $hg=g$. That is, $ G(f) \mu_D h= G(f) \mu_D$.
After we apply the functor $F$ and compose the morphism $\nu_M$ on the left,
we get
$$ \nu_M FG(f) F(\mu_D) F(h) =\nu_M FG(f) F(\mu_D).$$
The naturality of $\nu$ implies that $f \nu_{F(D)}=\nu_MFG(f)$. and so,
$ f \nu_{F(D)} F(\mu_D) F(h) =f \nu_{F(D)}F(\mu_D).$
On the other hand, we know that  $\nu_{F(D)} F(\mu_D)=1_{F(D)}$. Thus, $fF(h) =f$. As $f$ is an $F(D)$-cover, $F(h)$ is an isomorphism, and so is $h$ as $\mu_Dh=GF(h)\mu_D$ by the naturality of $\mu$. 
\cqfd

  \section{Relative PGF modules}\label{Sec. 4}

In this section, we introduce a new class of modules and investigate their main homological properties, aiming to prove Theorem A from the Introduction. 

\begin{defn}\label{PGCF def} Let $C$ be an $R$-module. An $R$-module $M$ is said to be \textbf{projectively coresolved $\GC$-flat}, or $\PGCF$ for short, if there exists a $\left(\Prod_R(C^+)\otimes_R-\right)$-exact exact complex of $R$-modules  $$\textbf{X}:\cdots \rightarrow P_1\rightarrow P_0 \rightarrow C_{-1}\rightarrow C_{-2}\rightarrow\cdots$$ with all $P_i$ projective,  $C_j\in\Add_R(C)$, and such that $M\cong \Im(P_0\to C_{-1})$.

The class of all projectively coresolved $\GC$-flat $R$-modules will be denoted by $\PGCF(R).$
\end{defn} 

\begin{rem}If $C=R$, or generally, $C$ is a projective generator $R$-module, then the class $\PGCF(R)$ coincides with ${\rm PGF}(R)$, the class of all {\rm PGF} modules recently introduced by \v{S}aroch and \v{S}\v{t}ov\'{\i}\v{c}ek in \cite{SS20}. 
\end{rem}

We begin our investigation with the following two propositions which immediately follow from the definition.
\begin{prop}\label{prop of PGcF 1}  The class $\PGCF(R)$ is closed under direct sums.
\end{prop}
\begin{prop}\label{prop of PGcF 2} An $R$-module $M$ is $\PGCF$ if and only if the following two assertions hold:
\begin{enumerate}
	\item[$(a)$] $\Tor_{i\geq 1}^R(\Prod_R(C^+),M)=0$.
	\item[$(b)$] $M$ has a $\left(\Prod_R(C^+)\otimes_R-\right)$-exact $\Add_R(C)$-coresolution
	$0 \rightarrow M\to  C_{-1}\rightarrow C_{-2}\rightarrow\cdots.$
\end{enumerate}

\end{prop}

\bigskip

Unlike the (absolute) Gorenstein modules, it is not immediately clear from the definition whether kernels in the exact sequence defining the relative Gorenstein modules are of the same type. However, in our setting, we have the following.

\begin{lem}\label{ker are PGcF} Let $C$ be $\Tor$-$\prod$-orthogonal and consider a $\left(\Prod_R(C^+)\otimes_R-\right)$-exact exact sequence $R$-modules $$\textbf{X}:\cdots \rightarrow P_1\rightarrow P_0 \rightarrow C_{-1}\rightarrow C_{-2}\rightarrow\cdots$$
with all $P_i$ projective and $C_j\in\Add_R(C)$. Then each kernel $K_j:=\Ker(C_{j}\rightarrow C_{j-1})\in\PGCF(R)$.

\end{lem}

\proof We only show $K_{-2}\in\PGCF(R)$ as the result will follow by induction.

By definition, $K_{-1}\in\PGCF(R)$ and so $\Tor^R_i(X,K_{-1})=0$ (Proposition \ref{prop of PGcF 2}) for every $X\in\Prod_R(C^+)$ and $i\geq 1$. Hence, 
$$\Tor^R_i(X,K_{-2})\cong \Tor^R_{i-1}(X,K_{-1})=0, \,\,\, \forall i\geq 2,$$
and so $\Tor^R_{i\geq 1}(X,K_{-2})=0$ as $\textbf{X}$ is a $\left(\Prod_R(C^+)\otimes_R-\right)$-exact exact sequence. On the other hand, it is clear that $K_{-2}$ has a $\left(\Prod_R(C^+)\otimes_R-\right)$-exact $\Add_R(C)$-coresolution
	$0 \rightarrow K_{-2}\to  C_{-2}\rightarrow C_{-3}\rightarrow\cdots.$
Therefore, applying Propsoition \ref{prop of PGcF 2} again, we obtain $K_{-2}\in\PGCF(R)$. \cqfd

\begin{prop}\label{Add-periodicity} The following two assertions hold.
	\begin{enumerate}
		\item[$(1)$] $C$ is $\Tor$-$\prod$-orthogonal if and only if  $\Add_R(C)\subseteq \PGCF(R)$.
		\item[$(2)$] If $C$ is $\Tor$-$\prod$-orthogonal, the following are equivalent for an $R$-module $M$:
		\begin{enumerate}
\item[$(a)$] $M\in\PGCF(R)$.
\item[$(b)$] There exists a short exact sequence $0\to M\to A \to G\to 0$ with $A\in\Add_R(C)$ and $G\in\PGCF(R)$. 
		\end{enumerate}
	
In other words, $\PGCF(R)$ is right $\Add_R(C)$-periodic.
	\end{enumerate}
\end{prop}
\proof $(1)$ The if-part is clear as in this case $C\in \Add_R(C)\subseteq \PGCF(R)$. Conversely, if $X\in\Add_R(C)$, then the exact complex $\cdots \to 0\to X\stackrel{1_X}{\to} X\to 0\to \cdots $ 
is a $\left(\Prod_R(C^+)\otimes_R-\right)$-exact exact complex of $R$-modules with $X$ in degrees $-1$ and $-2$. Therefore, $X\in\PGCF(R)$ according to Lemma \ref{ker are PGcF}.

(2) $(a)\Rightarrow (b)$ Follows by the definition and Lemma \ref{ker are PGcF}.

$(b)\Rightarrow (a)$ By Propsotion \ref{prop of PGcF 2} and using the short exact sequence $0\to M\to A \to G\to 0,$
we get $\Tor^R_i(X,M)\cong \Tor^R_{i+1}(X,G)=0$ for every $X\in\Prod_R(C^+)$ and every $i\geq 1$. Furthermore, the $(\Prod_R(C^+)\otimes_R-)$-exact $\Add_R(C)$-coresolution of $M$ can be constructed by combining the short exact sequence and a $(\Prod_R(C^+)\otimes_R-)$-exact $\Add_R(C)$-coresolution of $G$, which exists by Proposition \ref{prop of PGcF 2}. Hence, $M\in\PGCF(R)$ using Proposition \ref{prop of PGcF 2}(2) again.  \cqfd

\bigskip

One of the problems one encounters when dealing with Gorenstein objects that are defined by means of the tensor product is that it is not immediately clear whether the class of all these objects is closed under extensions. This closure property is found to be challenging for many classes of objects in Gorenstein homological algebra. Fortunately, we overcome this problem using an approach that characterizes objects defined via the tensor product using the Hom functor. This approach has been used successfully to settle the closure under extensions of the class of Gorenstein AC-projective modules (\cite[Theorem A.6]{BGH14}) and the class of PGF modules (\cite[Corollary 4.5]{SS20}).

In our situation, a key result, along with this approach, is a recent result of Gillespie and Iacob  \cite[Theorem 6]{GI22}, which is based on symmetric duality pairs. Given a class $\mathcal{L}$ of $R$-modules and a class $\mathcal{A}$ of right $R$-modules, $(\mathcal{L},\mathcal{A})$ is called a symmetric duality pair \cite[Definition 5]{GI22} if the follwing two assertions hold:
\begin{enumerate}
    \item $\forall \, _RM,$ $M\in \mathcal{L}$ if and only if $M^+\in \mathcal{A}$.
    \item $\forall N_R,$ $N\in \mathcal{A}$ if and only if $N^+\in \mathcal{L}$.
    \item $\mathcal{L}$ and $\mathcal{A}$ are closed under direct summands and finite direct sums.
\end{enumerate}

Unlike the pair $(\langle R\rangle,\langle R^+\rangle)$ (see \cite[Remark 2.1]{EIP20}), it is unkown whether $(\langle C\rangle,\langle C^+\rangle)$ is a symmetric duality pair. However, we still take advantage of \cite[Theorem 6]{GI22} to show the following.

\begin{prop}\label{new charac of PGcF}Assume that $C$ is pure projective and let $M$ be an $R$-module. The following assertions are equivalent.
	\begin{enumerate}
\item  $M\in \PGCF(R)$.
\item  There exists a $\Hom_R(-,\langle C\rangle)$-exact  exact sequence of $R$-modules  $$\textbf{X}:\cdots \rightarrow P_1\rightarrow P_0 \rightarrow C_{-1}\rightarrow C_{-2}\rightarrow\cdots$$  with all $P_i$ projective, $C_j\in\Add_R(C)$ and such that $M\cong \Im(P_0\to C_{-1})$.
\item The following two assertions hold: \begin{enumerate}
	\item[$(a)$] $\Ext^{i\geq 1}_R(M,A)=0$ for any $A\in \langle C\rangle$.
	\item[$(b)$]  $M$ has a $\Hom_R(-,\langle C\rangle)$-exact $\Add_R(C)$-coresolution; $0 \rightarrow M\to  C_{-1}\rightarrow C_{-2}\rightarrow\cdots.$
\end{enumerate}
	\end{enumerate}

In symbols, we have the following equalities $$\PGCF(R)=\GP_{ \langle C\rangle}(\proj(R),\Add_R(C))=\,^{\perp_\infty}\langle C\rangle \cap \GP_{ \langle C\rangle}(-,\Add_R(C)).$$
\end{prop}
\proof  $(2)\Leftrightarrow (3)$ Clear.

$(1) \Rightarrow (2)$ Assume that $M\in \PGCF(R)$. By definition, there exists a  $\left(\Prod_R(C^+)\otimes_R-\right)$-exact exact  complex of  $R$-modules $$\textbf{X}:\cdots \rightarrow P_1\rightarrow P_0 \rightarrow C_{-1}\rightarrow C_{-2}\rightarrow\cdots$$
with all $P_i$ projective and $C_j\in\Add_R(C)$ and such that $M\cong \Im(P_0\to C_{-1})$. Since $\langle C^+\rangle =\langle\Prod_R(C^+)\rangle$ (Lemma \ref{definability of Add and Prod}), it follows from Lemma \ref{def-tensor}(2) that $\XX$ is $(\langle C^+\rangle\otimes-)$-exact. Using now \cite[Theorem 1.2]{EIP20}, we get a symmetric duality pair $(\A,\langle C^+\rangle)$ with $\mathcal{A}=\Prod_R(\langle C^+\rangle^+)^p$. Then $\textbf{X}$ is $\Hom_R(-,\mathcal{A})$-exact by \cite[Theorem 6]{GI22}. Since $C^{++}\in \langle C^+\rangle^+$, it follows from Lemma \ref{definability of Add and Prod} that $\langle C\rangle=\langle C^{++}\rangle\subseteq \Prod_R(\langle C^+\rangle^+)^p=\mathcal{A}$. Consequently, $\textbf{X}$ is $\Hom_R(-,\langle C\rangle)$-exact. 

$(2)\Rightarrow (1)$ Given an exact sequence $\XX$ as in (2), we only need to show that $(C^+)^I\otimes \XX$  is exact for any index set $I$.  Clearly, we have the following natural isomorphisms of complexes
$$((C^+)^I\otimes \XX)^+\cong ((C^{(I)})^+\otimes_R \XX)^+\cong \Hom_R(\XX,(C^{(I)})^{++}).$$
Since the definable class $\langle C\rangle$ is closed under direct sums, $C^{(I)}\in\langle C\rangle$ and therefore
$(C^{(I)})^{++}\in\langle C\rangle$ by \cite[Lemma 1.1]{EIP20}. So, by hypothesis, the complex $$((C^+)^I\otimes \XX)^+\cong \Hom_R(\XX,(C^{(I)})^{++})$$ is exact and so is $(C^+)^I\otimes \XX$ as $\mathbb{Q}/\mathbb{Z}$ is an injective cogenerator $\mathbb{Z}$-module. \cqfd

\begin{cor}\label{PGcF are GC-proj} Assume that  $C$ is pure projective. Then, $\PGCF(R)\subseteq \GCP(R)$. Consequently, if $C$ is $\Tor$-$\prod$-orthogonal, then it is $\Sigma$-self-orthogonal.
\end{cor}
\proof The first claim follows by Proposition \ref{new charac of PGcF}, since $\Add_R(C)\subseteq \langle C\rangle$. The last claim follows by the fact that $C$ is $\Tor$-$\prod$-orthogonal (resp., $\Sigma$-self-orthogonal) if and only if $_RR\in\PGCF(R)$ (resp., $_RR\in\GCP(R)$). \cqfd

\begin{rem}Assume that $R$ is commutative and $C$ is semidualizing. Recall \cite[Definition 1.1]{ZWL14} that an $R$-module $M$ is ${\rm D}_C$-projective if it belongs to the class ${\rm D_CP}(R):=\GP_{\F(R)}(\proj(R),\Add_R(C))$. 
Since $\Add_R(C)\subseteq \F(R)\subseteq  \langle C\rangle$, 
       $\PGCF(R)\subseteq {\rm D_CP}(R)\subseteq \GCP(R)$ by Proposition \ref{new charac of PGcF}.
However, whether these inclusions are equalities is an open question even in the absolute case (see for instance \cite{ElM23}).
\end{rem}

Since the $\PGCF$ modules are now characterized by means of the $\Hom$ functor, more closure properties can easily be obtained.

\begin{prop}\label{clos-extens-dirsum-kerepi} Assume that $C$ is pure projective. Then, the class $\PGCF(R)$ is closed under extensions. If moreover, $C$ is $\Tor$-$\prod$-orthogonal, then the  class $\PGCF(R)$ is also closed under kernels of epimorphisms and direct summands.
\end{prop}

\proof By Proposition \ref{new charac of PGcF},  $\PGCF(R)=\GP_{ \langle C\rangle}(\proj(R),\Add_R(C))$ and $\Ext^1_R(\PGCF(R),\Add_R(C))=0$, and  so $\PGCF(R)$ is closed under extensions by Lemma \ref{cl-ext}. Moreover, if $C$ is $\Tor$-$\prod$-orthogonal, then $\PGCF(R)$ is right $\Add_R(C)$-periodic by Proposition \ref{Add-periodicity}. Hence, $\PGCF(R)$ is closed under kernels of epimorphisms by Lemma \ref{cl-epikernels}  and direct summands by Lemma \ref{cl-dir-summnds}. \cqfd

The last ingredient needed to prove Theorem A is when the class $\PGCF(R)$ contains all projective modules. Since $\PGCF(R)$ is closed under direct sums (Proposition \ref{prop of PGcF 1}) and direct summands (Proposition \ref{clos-extens-dirsum-kerepi}), we only need to verify when the class $\PGCF(R)$ contains the regular module $_RR$. This is precisely the second condition in the following definition.

\begin{defn}\label{g-tilting} An $R$-module $C$  is said to be \textbf{generalized Wakamatsu  tilting}, or \textbf{g-tilting} for short, if it satisfies the following two properties:
	\begin{enumerate}
\item $C$ is $\Tor$-$\prod$-orthogonal.
\item $R$ has a $\left(\Prod_R(C^+)\otimes_R-\right)$-exact $\Add_R(C)$-coresolution
$0 \rightarrow R\to  C_{-1}\rightarrow C_{-2}\rightarrow\cdots.$
	\end{enumerate}
\end{defn}

It is clear from this definition that $C$ satisfies condition (2) if and only if $R\in\PGCF(R),$ and that $C$ is $g$-tilting if and only if $R,C\in\PGCF(R)$. This observation will be used in several places in this paper.

\bigskip

Before we give the proof of Theorem A, it is of great importance to clarify the relation between g-tilting and other tilting-like modules that arise in homological algebra and tilting theory.

Recall that a module $C$ is called Wakamatsu tilting \cite{Wak04} if $(1)$ $_RC$ admits a degreewise finite projective resolution, $(2)$ $\Ext^{i \geq 1}_R(C,C) = 0$, and (3) there exists a $\Hom_R(-,C)$-exact coresolution $0 \to  R \to  C_{-1} \to C_{-2} \to \cdots,$  where each $C_j$ is isomorphic to a direct summand of finite direct sums of copies of $C$.

Wakamatsu showed in \cite[Corollary 3.2]{Wak04} that an $(R,S)$-bimodule $C$ is semidualizing if and only if $_RC$ is tilting with $S=\End_R(C)$ if and only if $C_S$ is tilting with $R= \End_S(C)$. From the following, we see that the generalized Wakamatsu tilting modules properly generalize both concepts. Note that any free module of infinite rank is generalized Wakamatsu tilting, but
it is neither semidualizing nor Wakamatsu tilting.
 
\begin{prop}\label{semid is g-tilting} All Wakamatsu tilting modules are g-tilting. Consequently, any  semidualizing $(R,S)$-bimodule is g-tilting both as a left $R$-module and as a right $S$-module
\end{prop}
\proof The proof is similar to that of \cite[Proposition 4.3]{BEGO22}. \cqfd

Recall that a module $C$ is called tilting (\cite[Definition *4]{Baz04}/\cite[Definition 13.1]{GT12}) if $(1)$ $_RC$ has a finite projective dimension, $(2)$ $\Ext^{i\geq 1}_R(C,C^{(I)}) = 0$ for any set $I$, and (3) there exists an $\Add_R(C)$-coresolution $0 \to  R \to  C_{-1}  \to \cdots \to C_{-r}\to 0$.

Rump \cite[Corollary 2]{Rum21} has shown recently that all tilting modules are w-tilting. In the following result, we show that tilting modules are also g-tilting.

\begin{prop}\label{til is g-til} All tilting modules are g-tilting. 
\end{prop}
\proof  Assume that $C$ is a tilting module. It follows from the proof of \cite[Proposition 3.2]{YYZ20} that $\Tor_{i\geq 1}^R((C^+)^{I},C)=0$ for any set $I$. Moreover, by definition, $_RR$ has an $\Add_R(C)$-coresolution 
$$0 \to  R \to  C_{-1} \to  \cdots \to C_{-r}\to 0.$$
It remains to show that this sequence is $(\Prod_R(C^+)\otimes_R-)$-exact. Let $X\in \Prod_R(C^+)$. Then, $\Tor_{i\geq 1}^R(X,C_j)=0$ for all $j=-1,-2\cdots ,-r-1$. So, if we set $K_j=\Ker(C_j\to C_{j-1})$, we get 
$$\Tor^R_1(X,K_j)\cong \Tor^R_2(X,K_{j-1}) \cong \cdots \cong \Tor^R_{r+j}(X,K_{-r+1})\cong \Tor^R_{r+j+1}(X,C_{-r})=0.  \text{ \cqfd} $$

\bigskip

\begin{prop}\label{w-til is g-til is w+-til} Any g-tilting module is $w^+$-tilting and any pure projective g-tilting module is w-tilting.
\end{prop}
\proof The first assertion follows from the definition as $\Add_R(C)\subseteq \F(R)$. The second assertion follows by Proposition \ref{new charac of PGcF} and the fact that $\Add_R(C)\subseteq \langle \Add_R(C)\rangle= \langle C\rangle$ by Lemma \ref{definability of Add and Prod}.
\cqfd

\bigskip

The following diagram summarizes the relations among all tilting-like modules that are related to g-tilting modules.
 $$\xymatrixcolsep{7pc}
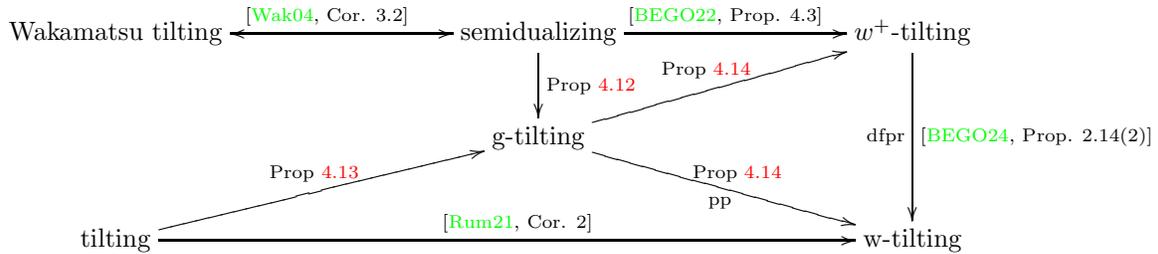
\xymatrix{  \text{Wakamatsu tilting} \ar[r]^{\text{\cite[Cor. 3.2]{Wak04}}} &  \text{semidualizing} \ar[l] \ar[d]^{\text{Prop \ref{semid is g-tilting} }}  \ar[r]^{\text{\cite[Prop. 4.3]{BEGO22}}} & \text{$w^+$-tilting}  \ar[dd]_{\text{dfpr}} ^{\text{\cite[Prop. 2.14(2)]{BEGO24}}}   \\    &   \text{g-tilting} \ar[ur]^{\text{Prop \ref{w-til is g-til is w+-til} }}  \ar[rd]_{\text{pp}}^{\text{Prop \ref{w-til is g-til is w+-til} }}    &   \\
  \text{tilting}  \ar[ru]^{\text{Prop \ref{til is g-til}}}  \ar[rr]^{\text{ \cite[Cor. 2]{Rum21}}}&        &   \text{w-tilting} } $$

\bigskip

Returning to the main result of this section, the following proposition is the final result needed to obtain Theorem A.

\begin{prop} \label{proj in PGcF}Let $C$ be pure projective. Then, $C$ is $g$-tilting if and only if $\proj(R)\cup\Add_R(C)\subseteq \PGCF(R)$.

\end{prop}
\proof  The "if" part is clear, since $R,C\in \proj(R)\cup\Add_R(C)\subseteq \PGCF(R)$. Conversely, if $C$ is g-tilting, then $R,C\in\PGCF(R)$. But, $\PGCF(R)$ is closed under direct sums (Proposition \ref{prop of PGcF 1}) and direct summands (Proposition \ref{clos-extens-dirsum-kerepi}). So, this desired inclusion follows. \cqfd

The following is Theorem A from the Introduction. 

\begin{thm}\label{PGCF is proj-resol} Let $C$ be a pure projective $R$-module. Then, $C$ is $g$-tilting if and only if $\PGCF(R)$ is projectively resolving and $\Add_R(C)\subseteq \PGCF(R)$. 
\end{thm}
\proof Follows by Propositions \ref{clos-extens-dirsum-kerepi} and  \ref{proj in PGcF}. \cqfd

We end this section with some useful consequences of this theorem. 

First, we obtain a new situation in which the class of $\GC$-flat modules is closed under extensions. Recall that whether the class of $\GC$-flat modules is closed under extensions is an open question, solved only under certain circumstances (see \cite{BEGO22,BEGO24}).
 
\begin{prop}\label{S is perfect} Assume that $C$ is finitely presented and $\prod$-$\Tor$-orthogonal. The following assertions are equivalent: 
	\begin{enumerate}
\item $S=\End_R(C)$ is left perfect. 
\item $\Add_R(C)=\F(R)$.
\item $\PGCF(R)=\GCF(R)$.
	\end{enumerate}
 In this case, the class of $\GC$-flat modules is closed under extensions.
\end{prop}

\proof  First of all, notice that $\F(R)=C\otimes_S\flat(S)$ and $\Add_R(C)=C\otimes_S\proj(S)$ by Lemma \ref{Fc-Pc-Ic}.

$(1)\Rightarrow (3) $ If $S$ is left perfect, then $\proj(S)=\flat(S)$ and hence $\Add_R(C)=\F(R)$. Thus,  $\PGCF(R)=\GCF(R)$ by Proposition \ref{prop of PGcF 2} and \cite[Proposition 4.8]{BEGO22}.

$(3)\Rightarrow (2)$ Let $X=C\otimes_S F\in\F(R)$. Then, $X\in \F(R)\subseteq \GCF(R)=\PGCF(R)$. Hence, $\Ext^1_R(X,\F(R))=0$ by Proposition \ref{new charac of PGcF} since $\F(R) \subseteq\langle C\rangle$. On the other hand, since $_SF$ is flat, there exists a pure short exact sequence of left $S$-modules
$$0\to K\to P\to F\to 0$$
with $P$ projective and $K$ flat. Therefore, the short exact sequence  $$0\to C\otimes_SK\to C\otimes_SP\to C\otimes_SF\to 0$$
splits as $C\otimes_S K\in \F(R)$. Therefore, $X=C\otimes_S F\in\Add_R(C)$.

$(2)\Rightarrow (1)$ Let $F$ be a flat left $S$-module. By hypothesis, $C\otimes_S F \in \F(R)=\Add_R(C)$. Hence $C\otimes_S F$ is a direct summand of some module $C^{(I)}$. On the other hand, we get by \cite[Theroem 3.2.14]{EJ00} that $F\cong S\otimes_S F\cong \Hom_R(C,C)\otimes_S F \cong \Hom_R(C,C\otimes_S F)$ since $_RC$ is finitely presented and $_SF$ is flat. Then $_SF$ is a direct summand of $\Hom_R(C,C^{(I)})\cong S^{(I)}$, hence projective.
  
Finally, the last claim follows by Theorem \ref{PGCF is proj-resol}. \cqfd

 Another application of Theorem \ref{PGCF is proj-resol}, is that iterating the procedure used to define the $\PGCF$ modules yields exactly the $\PGCF$ modules. Having this stability is crucial for proving Theorem B.   

%The first equality in the following proposition is a relative version of \cite[Theorem 5.7]{KS24} recently proved by Kaperonis and Stergiopoulou. However, our approach is different. We are mainly based on the work of Wang in \cite{Wan20} (see Lemma \ref{stability problem}).

 \begin{prop}\label{PGCF-stability}For any pure-projective  $g$-tilting module $C$, we have the equalities $$\PGCF(R)=\GF_{\Prod_R(C^+)}(\PGCF(R))=\GP_{\langle C\rangle}(\PGCF(R)).$$
\end{prop}
\proof  First of all note that $\Add_R(C)$ is a $\Prod_R(C^+)$-pure cogenerator for $\PGCF(R)$ by Proposition \ref{Add-periodicity}. On the other hand, since $\proj(R)\subseteq \PGCF(R)$ and $\PGCF(R)$ is closed under kernels of epimorphisms, $\proj(R)$ is a $\Prod_R(C^+)$-pure generator for $\PGCF(R)$. Therefore,  %So, Lemma \ref{stability problem} gives the the following equalities
%$$\GF_{\Prod_R(C^+)}(\proj(R),-)=\GF_{\Prod_R(C^+)}(\PGCF(R),-),$$  $$\GF_{\Prod_R(C^+)}(-,\Add_R(C))=\GF_{\Prod_R(C^+)}(-,\PGCF(R)).$$
\begin{eqnarray*}
\PGCF(R) &=&\GF_{\Prod_R(C^+)}(\proj(R),\Add_R(C)) \\
	&=&\GF_{\Prod_R(C^+)}(\proj(R),-)\cap \GF_{\Prod_R(C^+)}(-,\Add_R(C))\\
	&=& \GF_{\Prod_R(C^+)}(\PGCF(R),-)\cap  \GF_{\Prod_R(C^+)}(-,\PGCF(R)) \,\,\, \text{ (by Lemma \ref{stability problem})} \\
	&=& \GF_{\Prod_R(C^+)}(\PGCF(R),\PGCF(R)) \\
	&=& \GF_{\Prod_R(C^+)}(\PGCF(R)).
\end{eqnarray*} 

The equality $\PGCF(R)=\GP_{\langle C\rangle}(\PGCF(R))$ can be obtained in the same way using Lemma \ref{stability problem} together with the equality $\PGCF(R)=\GP_{\langle C\rangle}(\proj(R),\Add_R(C))$ from Proposition \ref{new charac of PGcF}. \cqfd

\section{Relative PGF dimension of modules.} \label{Sec. 5}

  Over commutative noetherian rings with respect to a semidualizing module $C$, Holm and J\o rgensen introduced the $\GC$-projective dimension in \cite{HJ06}. It has been considered the primary extension of the Foxby-Golod $\GC$-dimension to non-finitely generated modules. In this section we introduce the $\PGCF$ dimension and show that it offers an alternative definition for the $\GC$-projective dimension.

Given a class of $R$-modules, recall that an $R$-module $M$ is said to have $\X$-resolution dimension at most an integer $n \geq  0$, $\resdim_X(M)\leq  n$, if $M$ has a finite $\X$-resolution: 
$$0\to X_n\to \cdots \to X_1\to X_0\to M\to 0.$$
If $n$ is the least non-negative integer for which such a
sequence exists then its $\X$-resolution dimension is precisely $n$, and if there is no such $n$ then we define its $\X$-resolution dimension as $\infty$. We use $\overline{\X}^n$ to denote the class of all modules modules with $\X$-resolution dimension $\leq n$.
\begin{defn}
	The \textbf{$\PGCF$ dimension} of an  $R$-module $M$ is defined as the $\PGCF(R)$-resolution dimension of $M$:
            $$\PGCFpd_R(M):=\resdim_{\PGCF(R)}(M).$$

%We use \overline{\PGCF(R)} to denote the class of all $R$-modules with finite $\PGCF$  dimension
\end{defn}

 In addition to the projective dimension, there are two other homological dimensions that are of interest and are closely related to the $\PGCF$ dimension. Recall from \cite{BGO16} that the $\P$-projective and $\GC$-projective dimensions of an $R$-module $M$ are defined, respectively, as  
     $$\pcd_R(M):={\rm resdim}_{\Add_R(C)}(M) \text{ and }\GCpd(M)={\rm resdim}_{\GCP(R)}(M).$$

The proof of the following result follows by \cite[Theorem 1.1]{Hua22} along with Theroem  \ref{PGCF is proj-resol}.

\begin{prop} \label{PGCF dimesntion over ses}
	Let $C$ be pure projective g-tiling. Given a short exact sequence of $R$-modules $0\to K\to M\to N\to 0$, we have:
	\begin{itemize}
		\item[\rm (1)] If any two of $K$, $M$ or $N$ have finite $\PGCF$ dimension, then so has the third.
		\item[\rm (2)] $\PGCFpd_R(K)\leq \sup\{\PGCFpd_R(M),\PGCFpd_R(N)-1\},$ and the equality holds whenever $\PGCFpd_R(M)\neq \PGCFpd_R(N).$
		\item[\rm (3)] $\PGCFpd_R(M)\leq \sup\{\PGCFpd_R(K),\PGCFpd_R(N)\},$ and the equality holds whenever \\ $\PGCFpd_R(N)\neq \PGCFpd_R(K)+1.$
		\item[\rm (4)] $\PGCFpd_R(N)\leq \sup\{\PGCFpd_R(M),\PGCFpd_R(K)+1\},$ and the equality holds whenever  $\PGCFpd_R(K)\neq \PGCFpd_R(M).$
	\end{itemize}
\end{prop}

\begin{prop}\label{PGcF-pd =< n} Let $C$ be pure projective g-tiling. Given an integer $n\geq 0$, the followig are equivalent for an $R$-module $M$.
	\begin{enumerate}
\item  $\PGCFpd(M)\leq n.$
\item There exists a short exact sequence of $R$-modules 
$0\to M\to X\to G\to 0$ 
with $G\in\PGCF(R)$ and $\pcd_R(X)\leq n$. 
\item There exists a short exact sequence of $R$-modules 
 $0\to X\to G\to M\to 0$
with $G\in\PGCF(R)$ and $\pcd_R(X)\leq n-1$ $($if $n=0$, this is understood to mean $X=0$$)$.
	\end{enumerate}

Consequently, $\overline{\PGCF(R)}^{\,n}\cap \PGCF(R)^\perp=\overline{\Add_R(C)}^{\,n}$.

\end{prop}
\proof $(1) \Leftrightarrow (2)$ Similar to \cite[Propositions 3.3 and 3.4]{BGO16}.

$(1)\Rightarrow (3)$ The case $n=0$ is trivial, while the case $n\geq 1$ follows as in \cite[Theorem 3.5]{BGO16}.

$(3)\Rightarrow (1)$ By Proposition \ref{PGCF dimesntion over ses}(4) as $\PGCFpd_R(X)\leq \pcd_R(X)\leq n-1$ and $\PGCFpd_R(G)\leq n$.

It remains only to show the equality $\overline{\PGCF(R)}^{\,n}\cap \PGCF(R)^\perp=\overline{\Add_R(C)}^{\,n}$. 

$(\supseteq)$ Let $M\in\overline{\Add_R(C)}^{\,n}$. Then  $M\in  \overline{\PGCF(R)}^{\,n}$ as  $\Add_R(C)\subseteq \PGCF(R)$ by Proposition \ref{Add-periodicity}(1). On the other hand, consider a finite $\Add_R(C)$-resolution of $M$
$$0\to C_n\to C_{n-1}\to \cdots \to C_0\to M\to 0$$ 
and let $K_i=\Ker(C_i\to C_{i-1})$ with $C_{-1}=M$. Since $C_n\in\Add_R(C)\subseteq \PGCF(R)^\perp$, we get for any module $G\in\PGCF(R)$, the isomorphisms
$$\Ext_R^1(G,M)\cong \Ext_R^2(G,K_0)\cong \cdots \cong\Ext_R^{n}(G,K_{n-2})\cong \Ext_R^{n+1}(G,C_{n})=0.$$ 
Hence, $M\in\PGCF(R)^\perp$.

$(\subseteq)$ Conversely, assume $M\in \overline{\PGCF(R)}^{\,n}\cap \PGCF(R)^\perp$. By (3), we get a split short exact sequence
$0\to M\to X\to G\to 0$ with $X\in \overline{\Add_R(C)}^{\,n}$ and $G\in\PGCF(R)$. Thus, $M \in\overline{\Add_R(C)}^{\,n}$. \cqfd

\begin{rem} Note that the short exact sequence in Proposition \ref{PGcF-pd =< n}(3) is a special $\PGCF$-precover of $M$. Consequently, any module with finite $\PGCF$ dimension has a special $\PGCF$-precover when $C$ is pure projective g-tilting. 
\end{rem}
By the previous section (Proposition \ref{proj in PGcF} and Corollary \ref{PGcF are GC-proj}), we have the following containments
$$\proj(R)\cup\Add_R(C)\subseteq \PGCF(R)\subseteq \GCP(R).$$
Consequently, $\GCpd_R(M)\leq \PGCFpd_R(M)\leq min\{\pd_R(M),\pcd_R(M)\}$. The rest of this section is devoted to investigating how much these dimensions differ to each other.

In the next proposition, we show that the $\PGCF$ dimension is a refinement of the $\P$-projective dimension, while the $\GC$-projective dimension is a refinement of the $\PGCF$ dimension.

\begin{prop}\label{refinment} Let $C$ be pure projective g-tiling. 
\begin{enumerate}
\item If $M\in \PGCF(R)^\perp$, then $\PGCFpd_R(M)= \pcd_R(M)$. This is the case for any module $M$ with finite $\P$-projective dimension.
\item If $\PGCFpd_R(M)<\infty$, then $\GCpd(M)=\PGCFpd_R(M)$.
In particular, any $\GC$-projective with finite $\PGCF$ dimension is $\PGCF$.
\end{enumerate}
\end{prop}

\proof $(1)$ Given an integer $n\geq 0$, it follows from the equality in Proposition \ref{PGcF-pd =< n} that $\PGCFpd_R(M)\leq n$ if and only if $\pcd_R(M)\leq n$. Then, $\PGCFpd_R(M)= \pcd_R(M)$. For the last claim, note that if $M$ has finite $\P$-projective dimension, then $M\in \overline{\Add_R(R)}\subseteq\PGCF(R)^\perp$ by  Proposition \ref{PGcF-pd =< n}.

$(2)$ Clearly $\GCpd(M)\leq \PGCFpd_R(M)$ as $\PGCF(R)\subseteq \GCP(R)$ by Proposition \ref{new charac of PGcF}. For the other inequality, we may assume that $n=\GCpd_R(M)<\infty$. Then, there exists an exact sequence of $R$-modules
$$0\to G_n\to P_{n-1}\to \cdots \to P_0\to M\to 0$$
with each $P_i$ projective and $G_n\in \GCP(R)$. Since each $P_i\in\PGCF(R)$ by Proposition \ref{proj in PGcF}, we only need to show that $G_n\in\PGCF(R)$.  By hypothesis, $\PGCFpd(M)<\infty$. Then, $m:=\PGCFpd(G_n)<\infty$ using Proposition \ref{PGCF dimesntion over ses}. So, there exists by Proposition \ref{PGcF-pd =< n} a short exact sequence 
$\mathcal{E}:0\to X\to G\to G_n\to 0$ with $G\in\PGCF(R)$ and $\pcd_R(X)\leq m$. Using Proposition \ref{PGcF-pd =< n} again, we see that $X\in\overline{\Add_R(C)}^{\,m}\subseteq \PGCF(R)^\perp$, which implies that $\mathcal{E}$ is split and hence $G_n\in\PGCF(R)$ as desired.  \cqfd

White asked in \cite[Question 2.15]{Whi10} whether the $\GC$-projective dimension is a refinement of the projective dimension. Zhang, Zhu and Yan answered this question in \cite[Theorem 2.11]{ZZY13}, for the modules with degreewise finitely generated projective resolutions with respect to a faithful semidualizing $(R,S)$-bimodule $C$. Another partial answer was given by Bennis, Garcia Rozas and Oyonarte in \cite[Theorem 3.12(2)]{BGO16}. However, this is not true in general, as the following example shows.

\begin{exmp}$($\cite[Example 3.5]{BEGO23a} and \cite[ Example in pg. 63]{ARS97}$)$ \label{PGCF dim isn't a ref} Consider the quiver $$Q_n : v_1 \to v_2 \to v_3 \to  v_4$$  and let $R=kQ$ be the path algebra over an algebraic field $k$. For $i=1, \cdots, 4$, consider the only indecomposable injective modules over $R$:  $E_i:=\Hom_R(e_iR,k)$. Then, $C:= E_1\oplus E_2\oplus E_3\oplus E_4$ is a semidualizing $(R,R)$-bimodule by \cite[Example 3.5]{BEGO23a} and hence g-tilting by Proposition \ref{semid is g-tilting}. This implies that  $C\in \PGCF(R)$ and hence  $E_1\in \PGCF(R)$. However, $E_1$ has a finite projective dimension by \cite[Chap III, Prop 1.4]{ARS97} but is not projective as it corresponds to the representation
$$ k \to 0 \to  0 \to  0,$$
which is not projective. Thus, $\PGCFpd_R(E_1)=0\neq \pd_R(E_1).$

\end{exmp}

The following gives some necessary and sufficient conditions under which the $\GC$-projective (resp., $\PGCF$) dimension is a refinement of the projective dimension.

\begin{thm} \label{PGCF dim is a ref} Let $C$ be pure projective g-tiling. The following assertions are equivalent:
	\begin{enumerate}
\item $\PGCFpd_R(M)= \pd_R(M)$  for any $R$-module $M$ with finite projective dimension.
\item $\GCpd_R(M)= \pd_R(M)$ for any $R$-module $M$ with finite projective dimension.
\item Any $\PGCF$  module with finite projective dimension is projective. 
\item Any $\GC$-projective module with finite projective dimension is projective.
\item Any $\PGCF$ module with finite projective dimension is $\PGF$.
\item Any $\GC$-projective module with finite projective dimension is Gorenstein projective.
	\end{enumerate}
\end{thm}

\proof $(1)\Leftrightarrow (2)$ Let $M$ be an $R$-module such that $\pd_R(M)<\infty$. If (1) holds, then $\PGCFpd_R(M)<\infty$ and hence $\GCpd_R(M)=\PGCFpd_R(M)=\pd_R(M)$ by Proposition \ref{refinment}(2). The other direction follows from the inequality 
$\GCpd_R(M)\leq \PGCFpd_R(M)\leq \pd_R(M)$ as $\proj(R)\subseteq \PGCF(R)\subseteq \GCP(R)$ by Proposition \ref{proj in PGcF} and Corollary \ref{PGcF are GC-proj}. 

Now, we only prove the equivalences 
$(1)\Leftrightarrow (3) \Leftrightarrow (5)$
as $(2)\Leftrightarrow (4) \Leftrightarrow (6)$ can be obtained similarly.

$(1)\Rightarrow (3)$ If $M$ is a $\PGCF$ module with finite projective dimension, then $\pd_R(M)=\PGCFpd_R(M)=0$. That is, $M$ is projective.

$(3)\Rightarrow (5)$ Easy.

 $(5)\Rightarrow (1)$ Assume now that $\pd_R(M)<\infty$. The inequality $\PGCFpd_R(M)\leq \pd_R(M)$ is clear since $\proj(R)\subseteq \PGCF(R)$ by Proposition \ref{proj in PGcF}. To obtain the other inequality, we may assume that $n=\PGCFpd_R(M)<\infty.$ Then, there exists an exact sequence of $R$-modules
 $$0\to G_n\to P_{n-1}\to \cdots \to P_0\to M\to 0$$
with each $P_i$ projective and $G_n\in \PGCF(R)$. Since $\pd_R(M)<\infty$, $m:=\pd_R(G_n)<\infty$, i.e.,  $G_n$ has a finite projective resolution:
$$0\to Q_m\to Q_{m-1}\to \cdots \to Q_0\to G_n\to 0.$$
Let $K_i=\Ker(Q_i\to Q_{i-1})$ with $Q_{-1}=G_n$. Applying the functor $\Hom_R(G_n,-)$ to the short exact sequences 
$\mathcal{E}_i:0\to K_i\to Q_i\to K_{i-1}\to 0,$ 
we get 
$$\Ext_R^1(G_n,K_0)\cong \Ext_R^2(G_n,K_1)\cong \cdots \Ext_R^m(G_n,K_{m-1})= \Ext_R^m(G_n,Q_m)=0.$$
Since $Q_m\in \proj(R)\subseteq \PGF(R)^\perp$ and $G_n\in \PGF(R)$, the short exact sequence $\mathcal{E}_0:0\to K_0\to Q_0\to G_0\to 0$ splits and hence $G_n$ is projective.
\cqfd

\section{Homotopical aspects of $\PGCF$} %\label{Sec. 6}

In this section we prove Theorems A and B from the Introduction.  First, we need some key ingredients. We start with the closure of the class of $\PGCF$ modules under transfinite extensions.

\begin{lem}\label{PGCF-trans-exte} Assume that $C$ is $\Tor$-$\prod$-orthogonal pure projective. The class of projectively coresolved $\GC$-flat $R$-modules is closed under transfinite extensions.
\end{lem}
\proof  Note that $C$ is $\Sigma$-self-orthogonal by Corollary \ref{PGcF are GC-proj}. This implies that $\Add_R(C)$ is a self-orthogonal class and therefore closed under transfinite extensions by Lemma \ref{transf-ext 2}. In particular, this is the case for the class of projective modules $\proj(R)=\Add_R(R)$. So, since $$\PGCF(R)=\GF_{\Prod_R(C^+)}(\proj(R),-)\cap \GF_{\Prod_R(C^+)}(-,\Add_R(C)),$$ closure under transfinite extensions follows by Lemma \ref{transf-ext 1}.  \cqfd

For a chain complex $\XX=(X_n,d_n^\XX)$ and an integer $k$, let $\XX[k]$ denote the $k$-th shift of $\XX$. It is the chain complex whose degree $m$ is $\XX[k]_m=X_{m-k}$, and whose differentials are unchaged.

\begin{lem}\label{SPGCF} Let $M\in\PGCF(R)$. Then, $M$ is a direct summand in a module $N$ such that there exists a $\left(\Prod_R(C^+)\otimes_R-\right)$-exact exact sequence of $R$-modules $$0\to N\to P\oplus X\to N\to 0$$ with $P\in\proj(R)$ and $X\in \Add_R(C)$.

\end{lem}
\proof By definition, there exists an exact and  $\left(\Prod_R(C^+)\otimes_R-\right)$-exact complex 
     $$\XX:\cdots \rightarrow P_1\rightarrow P_0 \rightarrow C_{-1}\rightarrow C_{-2}\rightarrow\cdots$$ 
 with each $P_i\in\proj(R)$, $C_j\in\Add_R(C)$ and such that $M\cong \Im(P_0\to C_{-1})$. Adding all its shifts $\XX[k]$, we get a $(\Prod_R(C^+)\otimes_R-)$-exact exact  complex 
 $$\bigoplus_{k\in\mathbb{Z}}\XX[k]:\cdots\to P\oplus X\stackrel{f}{\to} P\oplus X\stackrel{f}{\to} P\oplus X \to\cdots $$ with $X:=\bigoplus _{j\leq -1} C_j\in \Add_R(C)$, $P:=\bigoplus_{i\geq 0} P_i \in \mathcal{P}(R)$ and such that $M$ is a direct summand of $N:=\Im(f)$. In this case, we have a $(\Prod_R(C^+)\otimes_R-)$-exact short exact sequence $0\to N\stackrel{u}{\to} P\oplus X\to N\to 0,$ as desired.  \cqfd

\bigskip

Now we are ready to prove Theorem B. The proof is inspired by that of \v{S}aroch and  \v{S}\v{t}ov\'{\i}\v{c}ek \cite{SS20}. For this result, we need a slightly stronger condition on our module $C$ than being pure projective. Namely, $C$ needs to be of type $FP_2$, i.e., $C$ and its first syzygy are finitely presented. Note that any module of type $FP_2$ is finitely presented and hence pure projective. The reason behind this assumption is that, in this case, the class $C^\perp$ is definable \cite[Example 6.10]{GT12}. This property is crucial in the proof of the following.

\begin{thm} \label{PGCF-cot-pair} Assume that $C$ is of type $FP_2$. Then, $C$ is g-tilting if and only if $(\PGCF(R),\PGCF(R)^\perp)$ is a hereditary complete cotorsion pair cogenerated by a set $\mathcal{S}\cup \{C\}$.
\end{thm}
\proof $(\Leftarrow)$ $C$ is g-tilting as $R\in{}^\perp(\PGCF(R)^\perp)= \PGCF(R) $ and $C\in {}^\perp((\mathcal{S}\cup \{C\})^\perp)=\PGCF(R)$

$(\Rightarrow)$ By \cite[Lemma 6.31]{GT12}, it is possible to find an infinite cardinal$\nu$ such that $R$ is a left $\nu$-noetherian ring, that is, each left ideal of $R$ is $\leq\nu-$generated.

Now, let $(\mathcal{A},\mathcal{B})$ be the cotorsion pair cogenerated by a representative set $\mathcal{S}$ of $\leq\nu$-presented ($=\,\leq\nu$-generated) modules from $\PGCF(R)$, that is, $(\mathcal{A},\mathcal{B})=(^\perp(\mathcal{S}^\perp),\mathcal{S}^\perp)$. By hypothesis, $R\in\mathcal{S}$. Then, $\mathcal{A}$ consists of all direct summands of $\mathcal{S}$-filtered modules by \cite[Corollary 6.14]{GT12}.  We claim that $(\PGCF(R),\PGCF(R)^\perp)$ is cogenerated by $\mathcal{S}$.  To see this, it suffices to show that $\PGCF(R)=\mathcal{A}$.

 $\bullet$ $\mathcal{A}\subseteq \PGCF(R)$: Since $\PGCF(R)$ is closed under direct summands by Proposition \ref{clos-extens-dirsum-kerepi} and transfinite extensions by Lemma \ref{PGCF-trans-exte}, we conclude that $\mathcal{A}\subseteq \PGCF(R)$.

 $\bullet$ $\PGCF(R)\subseteq \mathcal{A}$: Let $M\in\PGCF(R)$. By Lemma \ref{SPGCF}, $M$ is a direct summand in a module $N$ such that there exists a $\left(\Prod_R(C^+)\otimes_R-\right)$-exact exact sequence of $R$-modules $$\mathcal{E}:0\to N\stackrel{u}{\to} L\to N\to 0$$ with $L=P\oplus X$ such that $P\in\proj(R)$ and $X\in \Add_R(C)$.  Set $I=C^+$ and $\mathcal{D}=C^\perp$. Clearly, $\mathcal{D}$ contains the injective cogenerator $R^+$ of $R$-Mod and, by \cite[Example 6.10]{GT12}, is a definable class. It also contains the module $I^+$ since $\Ext^1_R(C,I^+)=\Ext^1_R(C,C^{++})\cong \Tor^R_1(C^+,C)^+=0$. Therefore, we obtain by \cite[Proposition 4.7(1)]{SS20} a filtration 
$$\mathfrak{F}=\{ \mathcal{E}_\alpha:0\to N_\alpha\stackrel{u_\alpha}{\to} L_\alpha\to N_\alpha\to 0\,|\;\alpha\leq \sigma \}$$ 
of $\mathcal{E}$ such that $L_{\alpha+1}/L_\alpha\in \,^\perp\mathcal{D}=\,^\perp(C^\perp)$ and $N_{\alpha+1}/N_\alpha$ is $\leq \nu $-presented  with   $N_{\alpha+1}/N_\alpha\in\,^\perp\langle I^+\rangle=\,^\perp\langle C^{++}\rangle=\,^\perp\langle C\rangle$ for every $\alpha<\sigma$. Note that the last equality follows by Lemma \ref{definability of Add and Prod}. It remains to show that $N_{\alpha+1}/N_\alpha\in\PGCF(R)$ for any $\alpha<\sigma$. Assembling the quotient short exact sequence
$$\mathcal{E}_{\alpha,\alpha+1}:=\mathcal{E}_{\alpha+1}/\mathcal{E}_{\alpha}\;:\; 0 \to N_{\alpha+1}/N_\alpha \to L_{\alpha+1}/L_\alpha \to  N_{\alpha+1} /N_\alpha \to  0$$ 
with itself, we get a $\Hom_R(-,\langle C\rangle)$-exact exact sequence of the form 
$$\XX\;:\;\cdots \to L_{\alpha+1}/L_\alpha\to L_{\alpha+1}/L_\alpha\to L_{\alpha+1}/L_\alpha\to \cdots .$$
On the other hand, since $_RC$ is a finitely presented module from $\PGCF(R)$, $C\in\mathcal{S}$ and hence $L_{\alpha+1}/L_\alpha\in \;^\perp(C^\perp)\subseteq \;^\perp(\mathcal{S}^\perp)=\mathcal{A}\subseteq \PGCF(R)$   for each $\alpha<\sigma$. Then, $N_{\alpha+1} /N_\alpha=\Im(L_{\alpha+1}/L_\alpha\to L_{\alpha+1}/L_\alpha)\in\GP_{\langle C\rangle}(\PGCF(R))=\PGCF(R)$ by Proposition \ref{PGCF-stability}.

Finally, the cotorsion pair $(\PGCF(R),\PGCF(R)^\perp)$ is hereditary as $\PGCF(R)$ is closed under kernels of epimorphisms by Proposition \ref{clos-extens-dirsum-kerepi} and the fact that it is complete follows by \cite[Theorem 6.11(b)]{GT12} as it is cogenerated by a set. \cqfd

 \begin{lem} \label{Add is covering} Assume that $C$ is finitely  presented. The following assertions are equivalent. 
 \begin{enumerate}
     \item[(1)]  $S=\End_R(C)$ is left perfect. 
     \item[(2)]  $\Add_R(C)$ is covering. That is, every $R$-module has a $\P$-projective cover.
     \item[(3)] Every $\F$-flat module has a $\P$-projective cover.
 \end{enumerate}
\end{lem}
\proof First of all consider the adjoint pair $(F,G)=(C\otimes_S -,\Hom_R(C,-)$. Since $C$ is finitely presented, it follows by \cite[Theroem 3.2.14]{EJ00} that the morphism $\mu_F: F\to \Hom_R(C,C\otimes_S F)$ is an isomorphism for any flat left $S$-module, in particular, for any module in $D:=\proj(S)$, the class of all left $S$-modules..

$(1)\Rightarrow (2)$ Since $S$ is left perfect, $\proj(S)=\flat(S)$ is covering (\cite[Theorem 5.3.1]{EJ00}). So, for any $R$-module $M$, the left $S$-module $G(M)=\Hom_R(C,M)$ has a projective cover. So, applying Lemma \ref{cov-env transfer}(1), we obtain that $M$ has an $\Add_R(C)$-cover, as 
$\Add_R(C)=C\otimes_S\proj(S)$  by Lemma \ref{Fc-Pc-Ic}(1). 

$(2)\Rightarrow (3)$ Clear.

$(3)\Rightarrow (1)$ By the proof of \cite[Theorem 5.3.2 $(1)\Leftrightarrow (2)$]{EJ00}, it suffices to show that every flat left $S$-module $F$ has a projective cover. By hypothsis,  $C\otimes_S F$ has an $\Add_R(C)$-cover. So,  applying again Lemma \ref{cov-env transfer}(1), we get that $F\cong \Hom_R(C,C\otimes_S F)=G(C\otimes_S F)$ has a projective cover.
\cqfd
\begin{cor}\label{PGCF is covering} Assume that $C$ is a g-tilting module of type $FP_2$. 
\begin{enumerate}
    \item $\PGCF(R)$ is a special precovering.
    \item $\PGCF(R)$ is covering if and only if  $S=\End_R(C)$ is left perfect.
\end{enumerate}
\end{cor}
\proof (1) Follows by Theorem \ref{PGCF-cot-pair}.

(2) $(\Rightarrow)$ By Lemma \ref{Add is covering}, it suffices to show that every module $C\otimes_S F\in \F(R)$ has an $\P$-projective cover. By hypothesis, the cotorsion pair $(\PGCF(R),\PGCF(R)^\perp)$ is perfect, and then $\PGCF(R)$ is covering.  
Let $F$ be a flat left $S$-module and $f: X\to C\otimes_S F$ be a $\PGCF$ cover. Let us show that $f$ is a $\P$-projective cover. Since $\PGCF(R)$ contains the projectives, $f$ is an epimorphism. So, by Wakamatsu Lemma (\cite[Lemma 5.13]{GT12}), $K:= \Ker f \in \PGCF(R)^\perp$. But $C\otimes_S F \in \F(R) \subseteq \langle C\rangle \subseteq \PGCF(R)^\perp$ by Proposition \ref{new charac of PGcF} and $\PGCF(R) \cap \PGCF(R)^\perp = \Add_R(C)$ by Proposition \ref{PGcF-pd =< n}. Hence, $X\in \PGCF(R) \cap \PGCF(R)^\perp = \Add_R(C)$. It remains to show that for any morphism $g :X\to X $ such that $fg = f$, is an automorphism. But this is the case as $f$ is a $\PGCF$ cover.

$(\Leftarrow)$ Conversely, if $S$
 is left perfect, then $\PGCF(R)=\GCF(R)$ by Corollary \ref{S is perfect}. Hence, $\PGCF(R)$ is closed under direct limits by \cite[Proposition 4.15]{BEGO22}. Thus, the cotorsion pair $(\PGCF(R),\PGCF(R)^\perp)$ is perfect by \cite[Corollary 5.32]{GT12}. Therefore, $\PGCF(R)$ is covering. \cqfd
\bigskip

One of the motivations behind the introduction of the class of $\PGCF$ modules is the lack of an answer to the following two natural questions, which are still open even in the absolute setting.

\begin{ques} Is any $\GC$-projective  module $\GC$-flat?
\end{ques} 

\begin{ques} Is the class of $\GC$-projective modules a (special) precovering class?
\end{ques}

As a consequence of Theorem \ref{PGCF-cot-pair}, we have partial answers to these two questions. 
\begin{cor} Assume that $C$ is g-tilting of type $FP_2$. Then 
 $\GCP(R)\subseteq \GCF(R)$
	if and only if $\GCP(R)=\PGCF(R)$. In this case, $(\GCP(R),\GCP(R)^\perp)$ is a complete and hereditary cotorsion pair, and consequently $\GCP(R)$ is a special precovering.

\end{cor}
\proof  $(\Rightarrow)$ The inclusion $\PGCF(R)\subseteq \GCP(R)$ follows by Corollary \ref{PGcF are GC-proj}. For the other inclusion, let $M\in\GCP(R)$.  Then, there exists a $\Hom_R(-,\Add_R(C))$-exact exact sequence 
$$\textbf{X}:\cdots \rightarrow P_1\rightarrow P_0 \rightarrow C_{-1}\rightarrow C_{-2}\rightarrow\cdots$$ 
 with each $C_j\in\Add_R(C)$ and $P_i\in\proj(R)$ and such that $M\cong \Im(P_0\to C_{-1})$. Since each image $I_i:=\Im(P_{i+1}\to P_i )\in \GCP(R)$ and each kernel $K_j:=\Ker(C_j\to C_{j+1})$ by \cite[Corollary 2.13]{BGO16}, $I_i,K_j\in \GCF(R)$ as $\GCP(R)\subseteq \GCF(R)$. Therefore, $\Tor_1^R(\Prod_R(C^+),I_i)=\Tor_1^R(\Prod_R(C^+),K_j)=0$. This implies that $\textbf{X}$ is $(\Prod_R(C^+)\otimes_R-)$-exact and hence $M\in\PGCF(R)$.

$(\Leftarrow)$ Follows from the fact that $\GCP(R)=\PGCF(R)\subseteq \GCF(R)$.

Finally, the last claim follows from Theorem \ref{PGCF-cot-pair}. \cqfd

Now we turn our attention to our last goal of this paper: Theorem C from the Introduction. Our main tool to construct the triple described in this theorem is a useful result by Gillespie \cite[Theorem 1.1]{Gil15} stating that if $(\Q,\widetilde{\R})$ and $(\widetilde{\Q},\R)$ are two complete hereditary cotorsion pairs in $\A$ such that  $\widetilde{\Q}\subseteq \Q$ and  $\Q\cap \widetilde{\R}=\widetilde{\Q}\cap \R$, then there exists a unique thick class $\W$ such that $(\Q,\W,\R)$ is a Hovey triple. We apply this result to the pairs: $$(\Q,\widetilde{\R})=(\PGCF(R),\PGCF(R)^\perp) \text{ and }(\widetilde{\Q},\R)=(^\perp\mathcal{B}_C(R),\mathcal{B}_C(R)).$$
We showed in Theorem \ref{PGCF-cot-pair} that the first pair is a complet and hereditary cotorsion pair. Now,  we will verify the remaining conditions.

\begin{prop}\label{Bass-cot-pair} Assume $C$ is  finitely presented. Then, $C$ is g-tilting if and only if  $(^\perp\mathcal{B}_C(R),\mathcal{B}_C(R))$
is a complete hereditary cotorsion pair cogenerated by a $\PGCF$ module. In this case, $^\perp\B_C(R)\subseteq \PGCF(R)$.
\end{prop}
\proof  $(\Leftarrow)$ If $(^\perp\mathcal{B}_C(R),\mathcal{B}_C(R))$ is a cotorsion pair cogenerated by a module $G\in\PGCF(R)$, then $\PGCF(R)^\perp\subseteq G^\perp:=\B_C(R)$ and hence $R,C\in {}^\perp\B_C(R)\subseteq \,^\perp(\PGCF(R)^\perp)= \PGCF(R)$. Then, $C$ is g-tilting.

$(\Rightarrow)$ By \cite[Theorem 3.10 and Corollary 3.12]{BDGO21}, the only thing we need to show is that the cotorsion pair $(^\perp\mathcal{B}_C(R),\mathcal{B}_C(R))$
is cogenerated by a $\PGCF$ module $G$.  Since $C$ is finitely presented g-tilting, it is w-tilting by Corollary \ref{PGcF are GC-proj}. Moreover, there exists a $(\Prod_R(C^+)\otimes_R-)$-exact $\Add_R(C)$-coresolution of $R$:
$$\textbf{X}\;:\;0\to R\stackrel{t_0}{\to} C_0\stackrel{t_1}{\to} C_1\stackrel{t_2}{\to} \cdots,$$
which is $\Hom_R(-,\Add_R(C))$-exact by Proposition \ref{new charac of PGcF} as $\Add_R(C)\subseteq\langle C\rangle$. Then, $\B_C(R)=H^{\perp_\infty}$ with $H=C\oplus(\oplus_{i\geq 0} \Coker t_i)$ by \cite[Theorem 3.10]{BDGO21}.  Note that $H\in\PGCF(R)$ as $C\in\PGCF(R)$ and each $\Coker t_i\cong \Ker t_{i+1}\in\PGCF(R)$ (Lemma \ref{ker are PGcF}).  On the other hand, consider a projective resolution of $H$ $$\cdots\to P_1\stackrel{f_1}{\to} P_0\stackrel{f_0}{\to} H\to 0$$ and let $K_{i+1}=\Ker(f_i)$. Since $\PGCF(R)$ contains all projective modules (Proposition \ref{proj in PGcF}) and is closed under kernels of epimorphisms (Proposition \ref{clos-extens-dirsum-kerepi}), each kernel $K_i\in\PGCF(R)$. Hence, $G=H\oplus (\oplus_{i\geq 1}K_i)\in\PGCF(R)$. Now, it is easy to see that, for any $i\geq 1$ and any $R$-module $A$, we have $\Ext^1_R(K_i,A)\cong \Ext^{i+1}_R(H,A)$. Then, 
$$\Ext^1_R(G,A)\cong \Ext^1_R(H,A)\oplus (\prod_{i\geq 1}\Ext^1_R(K_i,A))\cong\prod_{i\geq 1}\Ext^{i}_R(H,A).$$ 
Hence, $\B_C(R)=H^{\perp_\infty}=G^\perp$ with $G\in\PGCF(R)$ as desired.
\cqfd

\begin{prop}\label{core} Assume that $C$ is g-tilting and finitely presented. Then, 
	  $$\PGCF(R)\cap\PGCF(R)^\perp={}^\perp\mathcal{B}_C(R)\cap \mathcal{B}_C(R).$$
\end{prop}
\proof  Follows by \cite[Proposition 5.2]{BEGO24} 
and the equality in Proposition \ref{PGcF-pd =< n} for $n=0$. \cqfd

\begin{thm}\label{PGCF-Hov-trip}(\textbf{The $\PGCF$ model structure}) Assume that $C$ is of type $FP_2$. Then,  $C$ is g-tilting if and only if there exists a unique class $\W_C(R)$ such that $$\mathcal{M}_{\PGCF}=(\PGCF(R),\W_C(R),\mathcal{B}_C(R))$$ is a hereditary Hovey triple on $R$-Mod. In this case, there exists a cofibrantly generated hereditary abelian model structure on $R$-Mod such that:

\begin{enumerate}
\item[$(a)$] A morphism is a cofibration if and only if it is a monomorphism with cokernel in $\PGCF(R)$.

\item[$(b)$] A morphism is a fibration if and only if it is an epimorphism with kernel in $\B_C(R)$.

\item[$(c)$] A morphism is a weak equivalence if and only if it factors as a monomorphism with cokernel in $\W_C(R)$ followed by an
epimorphism with kernel kernel in $\W_C(R)$.
\end{enumerate}
\end{thm}
\proof $(\Leftarrow)$ $R\in {} ^\perp (\W_C(R)\cap \mathcal{B}_C(R))=\PGCF(R)$ and $C\in {}^\perp \B_C(R)=\PGCF(R)\cap \W_C(R)\subseteq \PGCF(R)$.

$(\Rightarrow)$ We know by Theorem \ref{PGCF-cot-pair} and Proposition \ref{Bass-cot-pair} that the pairs 
$$(\PGCF(R),\PGCF(R)^\perp) \text{ and } (^\perp\B_C(R),\B_C(R))$$
are complete and hereditary cotorsion pairs with ${}^\perp \B_C(R)\subseteq \PGCF(R)$. Furthermore, it follows by Proposition \ref{core} that these two cotorsion pairs have the same core, that is, 
$$\PGCF(R)\cap \PGCF(R)^\perp={}^\perp\B_C(R)\cap\B_C(R).$$
Therefore, following Gillespie \cite[Theorem 1.1]{Gil15}, there exists a unique thick class $\W_C(R)$ such that $(\PGCF,\W_C(R),\B_C(R))$ is a hereditary Hovey triple. 

Consequently, we obtain an abelian model structure, which is by \cite[Lemma 6.7 and Corollary 6.8]{Hov02}  cofibrantly generated in the sense of \cite[Sec 2.1.3]{Hov99} since the associated cotorsion pairs are each cogenerated by a set. Moreover, the descriptions of the cofibraions and fibrations follow by the definition, while the description of the weak equivalences follows by \cite[Lemma 5.8]{Hov02} and an application of the two out of three axiom (see also \cite[Proposition 2.3]{Gil16}).\cqfd

Assume that $C$ is g-tilting of type $FP_2$. One of the main consequences of Thoerem \ref{PGCF-Hov-trip} together with Gillespie \cite[Sections 5 and 6]{Gil11} (see also \cite[Proposition 4.2 and Theorem 4.3]{Gil16}), is that the full subcategory $\mathcal{A}_{c,f}:=\PGCF(R)\cap \mathcal{B}_C(R)$ is a Frobenius category (that is, an exact category with enough projectives and injectives, and in which the projective and injective objects coincide). In this case, the projective-injective objects are exactly the $\P$-projective modules. Thus, we have the  stable category 
$$\underline{\PGCF(R)\cap \mathcal{B}_C(R)}=(\PGCF(R)\cap \mathcal{B}_C(R))/\sim$$
where $\sim$ is the equivalence relation defined in $\PGCF(R)\cap \mathcal{B}_C(R)$ by $f \sim g$ if and only if  $f - g$ factors through a $\P$-projective module.  Furthermore, the homotopy category $\mathcal{M}_{\PGCF}$ of the $\PGCF$ model structure from Theorem~\ref{PGCF-Hov-trip} is triangle equivalent to the stable category $\underline{\PGCF(R)\cap \mathcal{B}_C(R)}$.

\bigskip

Unfortunately, \cite[Theorem 1.1]{Gil15} does not give an explicit description of the class of trivial objects. If $C=R$, or more generally, if $C$ is a projective generator, it follows by \v{S}aroch and \v{S}\v{t}ov\'{\i}\v{c}ek \cite[Theorem 4.9]{SS20}  that $\W_C(R)=\PGCF(R)^\perp$ since in this case $\PGCF(R)=\PGF(R)$. Following a recent development in \cite[Theorem D]{BEGO23b}, we obtain a partial answer under some finiteness of some global dimension of $R$.

Let us define the \textbf{global $\PGCF$ dimension} of $R$ as the supremum, if it exists, of the $\PGCF$ projective dimension of every $R$-module, that is,
$$\PGCFD(R):=\sup\{\PGCFpd_R(M)\,|\, \text{ $M$ is an $R$-module} \}.$$

Let us also define the \textbf{$^\perp\B_C$-projective dimension} of an $R$-module $M$ as the $^\perp\B_C(R)$-resolution dimension of $M$: $^\perp\B_C\,\text{-}\,\pd_R(M):=\resdim_{^\perp\B_C(R)}(M).$

The following lemma states that the $\PGCF$ dimension is a refinement of the $^\perp\B_C$-projective dimension. Note that if $M$ has finite $^\perp\B_C$-projective dimension, then $M\in\W_C(R)$ as $^\perp\B(R)\subseteq \W_C(R)$ and $\W_C(R)$ is thick.
\begin{lem}\label{PGCF dim is a ref 2}Assume that $C$ is g-tilting of type $FP_2$ and let $M$ be an $R$-module. Then, $\PGCFpd_R(M)\leq {}^\perp\B_C\,\text{-}\,\pd_R(M)$ with equality whenever $M\in\W_C(R)$.
\end{lem}
\proof The inquality follows since $^\perp\B_C(R)\subseteq \PGCF(R)$ by Proposition \ref{Bass-cot-pair}. Now, let $M\in\W_C(R)$. We may assume that $n:=\PGCFpd(M)<\infty$ and consider a partial projective resolution of $M$
$$0\to K_n\to P_{n-1}\to \cdots \to P_1 \to P_0\to M\to 0.$$
As each $P_i\in\proj(R)\subseteq \PGCF(R)$, it follows by Proposition \ref{PGCF dimesntion over ses}(2) that  $K_{n}\in\PGCF(R)$. But since  $P_i\in{}^\perp\B_C(R)\subseteq \W_C(R)$ and $\W_C(R)$ is thick, $K_n\in\PGCF(R)\cap \W_C(R)={}^\perp \B_C(R)$. Thus $^\perp\B_C\,\text{-}\,\pd_R(M)\leq n=\PGCFpd_R(M)$.  \cqfd
\begin{prop} \label{trivial objects} Assume $C$ is g-tilting of type $FP_2$. Then  $\PGCFD(R)\leq n$ if and only if the following two assertions hold. \begin{enumerate}
        \item[(a)] $(\PGCF(R), \overline{^\perp\B_C(R)}^n, \B_C(R))$ is a hereditary Hovey triple.
        \item[(b)] $\id_R(\langle C\rangle):=\sup\{ \id_R(M)| M\in \langle C\rangle\}\leq n$.
    \end{enumerate}

\end{prop}
\proof $(\Rightarrow)$ (a) By Theorem \ref{PGCF-Hov-trip} and \cite[Theorem D]{BEGO23b}, we obtain a hereditary Hovey triple $(\PGCF(R), \overline{^\perp\B_C(R)}, \B_C(R))$. Clearly, $\overline{^\perp\B_C(R)}^n\subseteq \overline{^\perp\B_C(R)}$. To see the other inclusion, let $M\in \overline{^\perp\B_C(R)}$. Then $^\perp\B_C-\pd_R(M)=\PGCFpd(M)\leq n$ by Lemma \ref{PGCF dim is a ref 2}.

Let us now show $(b)$. Let $M\in \langle C\rangle$ and $N$ be an $R$-module. By hypothesis, $\PGCFpd(N)\leq n$, that is, there exists a finite $\PGCF$ resolution of $N$:
$$0\to G_n\to \cdots \to G_0\to N\to 0.$$
By Proposition \ref{new charac of PGcF}, $\Ext_R^{i\geq 1}(G_k,M)=0$ for each $k$. So,
$\Ext_R^i(N,M)\cong \Ext_R^{i-n}(G_n,M)=0\, \forall i\geq n+1$. This means that $\id_R(M)\leq n$.

$( \Leftarrow)$ Consider a projective and an injective resolution of an $R$-module $M$: $$\cdots \to P_1\to P_0\to M\to 0 \text{ and } 0\to M\to I_0\to I_1\to\cdots,$$ 
respectively. Decomposing these exact sequences into short exact ones we get, for every integer $i\in \N$, 
	$$0 \to L_{i+1}\to P_i\to L_i\to 0\text{ \;\;and\;\;  }0\to K_i\to I_i\to K_{i+1}\to 0$$ 
	where $L_i=\Coker(P_{i+1}\to P_{i})$ and  $K_i=\Ker(I_i\to I_{i+1})$. Note that $M=L_0=K_0$. 
	Adding the direct sum of the first sequences, 
	$$0 \to \bigoplus_{i\in\N}L_{i+1}\to \bigoplus_{i\in\N}P_i\to M\oplus( \bigoplus_{i\in\N}L_{i+1})\to 0$$
	to the direct product of the second ones, 
	$$0\to  M\oplus( \prod_{i\in\N}K_{i+1})\to \prod_{i\in\N}I_i\to \prod_{i\in\N}K_{i+1}\to 0$$
	we get a short sequence of $R$-modules 
	$0\to N\to X\to N\to 0$, where $$X=(\bigoplus_{i\in N} P_i)\oplus (\prod_{i\in N} I_i)\text{ and } N=M\oplus \left((\bigoplus_{i\in N} L_{i+1}) \oplus (\prod_{i\in N} K_{i+1})\right).$$
Now, if $0\to D_n\to Q_{n-1} \to \cdots  \to Q_1\to Q_0\to N\to $ is a partial projective resolution of $N$, using  Horseshoe lemma, we construct a partial projective resolution of $X$: $$0\to G_n\to  Q_{n-1}\oplus Q_{n-1}\to   \cdots \to Q_1\oplus Q_1\to Q_0\oplus Q_0\to X\to  0$$ and a short exact sequence of $R$-modules:  $$0\to D_n\to G_n\to D_n\to 0.$$
Since $\prod_{i\in N} I_i$ is injective, $\prod_{i\in N} I_i\in\W_C(R)$, which implies by Lemma \ref{PGCF dim is a ref 2} that $\PGCFpd_R(\prod_{i\in N} I_i)={}^\perp\B_C-\pd_R(\prod_{i\in N} I_i)\leq n$ and therefore $\PGCFpd_R(X)\leq n$ as $\oplus_{i\in N} P_i\in \proj(R)\subseteq \PGCF(R)$. Since each $Q_i\oplus Q_i\in \PGCF(R)$, this implies by Proposition \ref{PGCF dimesntion over ses}(2) that $G_n\in \PGCF(R)$. On the other hand, since $\id_R(X)\leq n$ for all $X\in\langle C \rangle$ by (b), we get $\Ext_R^1(L_n,X)\cong \Ext_R^{1+n}(N,X)=0$. Therefore, gluing the short exact sequence $0\to D_n\to G_n\to D_n\to 0$ with itself infinitely many times, we get a $\Hom_R(-,\langle C \rangle)$-exact sequence of $R$-modules
$\cdots \to G_n\to G_n\to G_n\to \cdots.$ 
Therefore,  it follows by Proposition \ref{PGCF-stability} that $L_n\in \GP_{\langle C\rangle}(\PGCF(R))=\PGCF(R)$. Thus, $\PGCFpd(M)\leq \PGCFpd(N)\leq n$. \cqfd
\bigskip

We close this paper with some questions related to the dual situation. 

In the absolute case, one can define a dual concept of PGF modules (as in \cite[pg. 19]{ElM23}). However, it follows from Gillespie and Iacob \cite[Proposition 13(1)]{GI22} that these modules are nothing but Ding injective modules. 

In the relative setting, Ding injective modules with respect to a semidualizing module $C$ over commutative rings, $D_C$-injective, were introduced by Zhang, Wang and  Liu in \cite{ZWL14}. Their definition can be extended to the non-commutative setting with respect to a non-necessarily semidualizing module as follows.

\begin{defn}\label{DC-injective} Given an $R$-module $U$, an $R$-module $M$ is said to be \textbf{$D_U$-injective} if there exists a $\Hom_R({\rm Cogen_*}(U),-)$-exact exact sequence of $R$-modules 
$$\textbf{Y}:\;\cdots \to U_1\to U_0\to E_{-1}\to E_{-2}\to \cdots,$$ with each $U_i\in\Prod_R(U)$, $E_j\in\inj(R)$, and such that $M\cong \Im(U_0\to E_{-1})$.
\end{defn}

Note that we have seen by Lemma \ref{Fc-Pc-Ic}(2) that $\Prod_R(C^+)=\Hom_S(C,\mathcal{I}(S))$. On the other hand, using a similar type of arguments as in \cite[Proposition 3.3]{BEGO22}, it is straightforward to show that $\Cogen^*_R(C^+)=\Hom_S(C,\mathcal{FI}(S))$
where $\mathcal{FI}(S)$ denotes the class of all FP-injective right $S$-modules. Then, an $R$-module $M$ is ${\rm D_C}$-injective in the sense of \cite{ZWL14} if and only if it is ${\rm D_{U}}$-injective with $U=C^+$ in the sense of Definition \ref{DC-injective}.

 Now, looking at the ${\rm D_U}$-injective modules as a dual concept of the $\PGCF$ modules, a dual question to the one addressed in Theorem \ref{PGCF-Hov-trip} can be asked.

\begin{ques} Given an $R$-module $U$, when is ${\rm D_UI}(R)$ the class of all ${D_U}$-injective modules the class of cofibrant objects in a Hovey triple $(\Q,\W,{\rm D_UI}(R))?$
\end{ques}

In the absolute case, this question has recently been solved by Gillespie and Iacob \cite[Theorem 44]{GI22} while this problem seems to be more challenging in the relative setting.

Assume that $C$ is self-cosmall and $w^+$-tilting. It follows from \cite[Propositions 5.1(2) and 5.2(2)]{BEGO24}, that $(\A_C(R),\A_C(R)^\perp)$ is a complete and
hereditary cotorsion pair in Mod-$R$ with $\A_C(R)\cap \A_C(R)^\perp=\Prod_R(C^+)$. Furthermore, by a dual argument to Propositions \ref{Bass-cot-pair} and \ref{core}, one can show that $^\perp{\rm D_{C^+}I}(R)\cap {\rm D_{C^+}I}(R)=\Prod_R(C^+)=\A_C(R)\cap \A_C(R)^\perp$ and $ ^\perp{\rm D_{C^+}I}(R)\subseteq \A_C(R)$. Consequently, $(\Q,\W,{\rm D_{C^+}}I(R))$ is a Hovey triple for some thick class $\W$ if and only if $(^\perp {\rm D_{C^+}I}(R),{\rm D_{C^+}I}(R))$ is a complete cotorsion pair. So, the above question is reduced to:

\begin{ques}  Given an $R$-module $U$, when is ${\rm D_UI}(R)$ the right hand-side of a complete cotorsion pair $(^\perp {\rm D_UI}(R),{\rm D_UI}(R))$?
\end{ques}
\bigskip

\noindent\textbf{Acknowledgement.} The author is grateful to Jan \v{S}aroch  for providing proof of Lemma \ref{definability}.

\end{document}